\documentclass[11pt,english]{article}
\usepackage[T1]{fontenc}
\usepackage[latin1]{inputenc}
\usepackage{graphicx}

\makeatletter

\providecommand{\tabularnewline}{\\}

\usepackage{amsthm,amsfonts,amsmath,amssymb}

\theoremstyle{definition}
\newtheorem{defn}{Definition}
\theoremstyle{plain}
\newtheorem{thm}[defn]{Theorem}
\theoremstyle{plain}
\newtheorem{lem}[defn]{Lemma}
\theoremstyle{plain}
\newtheorem{cor}[defn]{Corollary}
\theoremstyle{definition}
\newtheorem*{rem*}{Remark}
\theoremstyle{plain}
\newtheorem{observation}[defn]{Observation}

\usepackage{fullpage,enumerate}

\newcommand{\dtv}{\mathrm{d}_\mathrm{TV}}

\def\epsilon{\varepsilon}

\def\Mscan{\mathcal{M}_{ \rightarrow }}

\def\Pscan{P_{\rightarrow}}

\newcommand{\Es}[2]{\textnormal{\textbf{E}}_{ #1 } \left[ #2 \right]}

\def\Prob{\textnormal{Pr}}

\newcommand{\R}{\mathbb{R}}

\def\Pk{P^{[k]}}
\def\Pj{P^{[j]}}
\def\PI{P^{[i]}}
\def\P1{P^{[1]}}
\def\Pm{P^{[m]}}

\def\M{\mathcal{M}}

\def\scan{\mathcal{M}_\rightarrow}

\def\escan{\mathcal{M}_\textnormal{edge}}

\def\treeblock{\mathcal{M}_\textnormal{tree}}

\def\Mix{\textnormal{Mix}}

\usepackage{babel}
\makeatother
\begin{document}

\title{Dobrushin Conditions for Systematic Scan with Block Dynamics%
\thanks{This work was partly funded by EPSRC projects GR/T07343/02 and GR/S76168/01. %
}}

\author{Kasper Pedersen\\
Department of Computer Science\\
University of Liverpool\\
Liverpool L69 3BX, UK\\
k.pedersen@csc.liv.ac.uk}
\maketitle
\begin{abstract}
We study the mixing time of systematic scan Markov chains on finite
spin systems. It is known that, in a single site setting, the mixing
time of systematic scan can be bounded in terms of the influences
sites have on each other. We generalise this technique for bounding
the mixing time of systematic scan to block dynamics, a setting in
which a (constant size) set of sites are updated simultaneously. In
particular we consider the parameter $\alpha$, corresponding to the
maximum influence on any site, and show that if $\alpha<1$ then the
corresponding systematic scan Markov chain mixes rapidly. As applications
of this method we prove $O(\log n)$ mixing of systematic scan (for
any scan order) for heat-bath updates of edges for proper $q$-colourings
of a general graph with maximum vertex-degree $\Delta$ when $q\geq2\Delta$.
We also apply the method to improve the number of colours required
in order to obtain mixing in $O(\log n)$ scans for systematic scan
for heat-bath updates on trees, using some suitable block updates. 
\end{abstract}

\section{Introduction\label{sec:Introduction}}

This paper is concerned with the study of finite \emph{spin systems}.
A spin system is composed of a set of sites and a set of spins, both
of which will be finite throughout this paper. The interconnection
between the sites is determined by an underlying graph. A configuration
of the spin system is an assignment of a spin to each site. If there
are $n$ sites and $q$ available spins then this gives rise to $q^{n}$
configurations of the system, however some configurations may be illegal.
The specification of the system determines how the spins interact
with each other at a local level, such that different local configurations
on a subset of the graph may have different relative likelihoods.
This interaction hence specifies a probability distribution, $\pi$,
on the set of configurations. One class of configurations that receive
much attention in theoretical computer science is \emph{proper $q$-colourings}
of graphs. A proper colouring is a configuration where no two adjacent
sites are assigned the same colour. One important example of a spin
system is when the set of legal configurations is the set of all proper
$q$-colourings of the underlying graph and $\pi$ is the uniform
distribution on this set. In statistical physics the spin system corresponding
to proper $q$-colourings is known as the $q$-state anti-ferromagnetic
Potts model at zero temperature.

Sampling from $\pi$ is a computationally challenging task. It is,
however, an important one and is often carried out by simulating some
suitable random \emph{dynamics} on the set of configurations. Such
a dynamics must have the following two properties

\begin{enumerate}
\item the dynamics eventually converges to $\pi$, and 
\item the rate of convergence (\emph{mixing time}) is polynomial in the
number of sites. 
\end{enumerate}
It is generally straightforward to ensure that a dynamics converges
to $\pi$ but much harder provide good upper bounds on the rate of
convergence, which is what we will be concerned with in this paper.

Arguably the simplest dynamics is the heat-bath Glauber dynamics which,
at each step, selects a site uniformly at random and updates the spin
assigned to that site by drawing a new spin from the distribution
on the spin of the selected site induced by $\pi$. This procedure
is repeated until the distribution of the Markov chain is sufficiently
close to $\pi$ using some suitable measure of closeness between probability
distributions. This dynamics falls under a family of Markov chains
that we call \emph{random update Markov chains}. We say that a Markov
chain is a random update Markov chain if the sites are updated in
a random order. This type of Markov chain has been frequently studied
in theoretical computer science and much is known about the mixing
time of various random update Markov chains.

An alternative to random update Markov chains is to construct a Markov
chain that cycles through (and updates the spin according to the induced
distribution) the sites (or subsets of sites) in a deterministic order.
We call this a \emph{systematic scan Markov chain} (or systematic
scan for short). Although systematic scan updates the sites in a deterministic
order it remains a random process since the procedure used to update
the spin assigned to a site is randomised, as specified by the appropriate
induced distribution. Systematic scan may be more intuitively appealing
that random update in terms of implementation, however until recently
little was know about the convergence rates of this type of dynamics.
It remains important to know how many steps one needs to simulate
a systematic scan for in order to for it to become sufficiently close
to its stationary distribution and recently there has been an interest
among computer scientists in investigating various approaches for
analysing the mixing time of systematic scan Markov chains, see e.g.
Dyer, Goldberg and Jerrum~\cite{dobrushin_scan,systematic_scan}
and Bordewich, Dyer and Karpinski~\cite{metric}. In this paper we
present a new method for analysing the mixing time of systematic scan
Markov chains, which is applicable to any spin system. As applications
of this method we improve the known parameters required for rapid
mixing of systematic scan on 

\begin{enumerate}
\item proper colourings of general graphs and 
\item proper colourings of trees. 
\end{enumerate}
A key ingredient in our method for proving mixing of systematic scan
is to work with a \emph{block dynamics}. A block dynamics is a dynamics
in which we allow a set of sites to be updated simultaneously as opposed
to updating one site at a time as in the description of the Glauber
dynamics above. Block dynamics is not a new concept and it was used
in the mid 1980s by Dobrushin and Shlosman \cite{DS} in their study
of conditions that imply uniqueness of the Gibbs measure of a spin
system, a topic closely related to studying the mixing time of Markov
chains (see for example Weitz's PhD thesis \cite{weitz_thesis}).
More recently, a block dynamics has been used by Weitz \cite{dror_combinatorial}
when, in a generalisation of the work of Dobrushin and Shlosman, studying
the relationship between various \emph{influence parameters} (also
in the context of Gibbs measures) within spin systems and using the
influence parameters to establish conditions that imply mixing. Using
an influence parameter to establish a condition which implies mixing
of systematic scan is a key aspect of the method presented in this
paper as we will discuss below. Dyer, Sinclair, Vigoda and Weitz \cite{mixing_time_space}
have also used a block dynamics in the context of analysing the mixing
time of a Markov chain for proper colourings of the square lattice.
Both of these papers consider a random update Markov chain, however
several of ideas and techniques carry over to systematic scan as we
shall see.

We will bound the mixing time of systematic scan by studying the \emph{influence}
that the sites of the graph have on each other. This technique is
well-known and the influence parameters generalised by Weitz \cite{dror_combinatorial}:
{}``the influence \emph{on} a site is small'' (originally attributed
to Dobrushin \cite{D}) and {}``the influence \emph{of} a site is
small'' (originally Dobrushin and Shlosman \cite{DS}) both imply
mixing of the corresponding random update Markov chain. It is worth
pointing out that a condition of the form ``if the influence \emph{on}
a site is small then the corresponding dynamics converges to $\pi$
quickly'' is known as a Dobrushin condition. In the context of systematic
scan, Dyer et al. \cite{dobrushin_scan} point out that, in a single
site setting, the condition {}``the influence \emph{on} a site is
small'' implies rapid mixing of systematic scan. Our method for proving
rapid mixing of systematic scan is a generalisation of this influence
parameter to block dynamics.

We now formalise the concepts above and state our results. Let $C=\{1,\dots,q\}$
be the set of spins and $G=(V,E)$ be the underlying graph of the
spin system where $V=\{1,\ldots,n\}$ is the set of sites. We associate
with each site $i\in V$ a positive weight $w_{i}$. Let $\Omega^{+}$
be the set of all configurations of the spin system and $\Omega\subseteq\Omega^{+}$
be the set of all legal configurations. Then let $\pi$ be a probability
distribution on $\Omega^{+}$ whose support is $\Omega$ i.e., $\{ x\in\Omega^{+}\mid\pi(x)>0\}=\Omega$.
If $x\in\Omega^{+}$ is a configuration and $j\in V$ is a site then
$x_{j}$ denotes the spin assigned to site $j$ in configuration $x$.
For each site $j\in V$, let $S_{j}$ denote the set of pairs $(x,y)\in\Omega^{+}\times\Omega^{+}$
of configurations that only differ on the spin assigned to site $j$,
that is $x_{i}=y_{i}$ for all $i\neq j$.

We will use Weitz's~\cite{dror_combinatorial} notation for block
dynamics, although we only consider a finite collection of blocks.
Define a collection of $m$ blocks $\Theta=\{\Theta_{k}\}_{k=1,\dots,m}$
such that each block $\Theta_{k}\subseteq V$ and $\Theta$ covers
$V$, where we say that $\Theta$ covers $V$ if $\bigcup_{k=1}^{m}\Theta_{k}=V$.
One site may be contained in several blocks and the size of each block
is not required to be the same, we do however require that the size
of each block is bounded independently of $n$. For any block $\Theta_{k}$
and a pair of configurations $x,y\in\Omega^{+}$ we write {}``$x=y$
on $\Theta_{k}$'' if $x_{i}=y_{i}$ for each $i\in\Theta_{k}$ and
similarly {}``$x=y$ off $\Theta_{k}$'' if $x_{i}=y_{i}$ for each
$i\in V\setminus\Theta_{k}$. We also let $\partial\Theta_{k}=\{ i\in V\setminus\Theta_{k}\mid\exists j\in\Theta_{k}:\{ i,j\}\in E(G)\}$
denote the set of sites adjacent to but not included in $\Theta_{k}$;
we will refer to $\partial\Theta_{k}$ as the \emph{boundary} of $\Theta_{k}$.

With each block $\Theta_{k}$, we associate a transition matrix $\Pk$
on state space $\Omega^{+}$ satisfying the following two properties: 

\begin{enumerate}
\item \label{item:update_theta}If $\Pk(x,y)>0$ then $x=y$ off $\Theta_{k}$,
and also 
\item \label{item:update_pi}$\pi$ is invariant with respect to $\Pk$. 
\end{enumerate}
Property \ref{item:update_theta} ensures that an application of $\Pk$
moves the state of the system from from one configuration to another
by only updating the sites contained in the block $\Theta_{k}$ and
Property \ref{item:update_pi} ensures that any dynamics composed
solely of transitions defined by $\Pk$ converges to $\pi$. While
the requirements of Property \ref{item:update_theta} are clear we
take a moment to discuss what we mean in Property \ref{item:update_pi}.
Consider the following two step process in which some configuration
$x$ is initially drawn from $\pi$ and then a configuration $y$
is drawn from $\Pk(x)$ where $\Pk(x)$ is the distribution on configurations
resulting from applying $\Pk$ to a configuration $x$. We than say
that $\pi$ is invariant with respect to $\Pk$ if for each configuration
$\sigma\in\Omega^{+}$ we have $\Pr(x=\sigma)=\Pr(y=\sigma)$. That
is the distribution on configurations generated by the two-step process
is the same as if only the first step was executed. In terms of our
dynamics this means that once the distribution of the dynamics reaches
$\pi$, $\pi$ will continue be the distribution of the dynamics even
after applying $\Pk$ to the state of the dynamics. Our main result
(Theorem \ref{thm:main_d}) holds for any choice of update rule $\Pk$
provided that it satisfies these two properties.

The distribution $\Pk(x)$, which specifies how the dynamics updates
block $\Theta_{k}$, clearly depends on the specific update rule implemented
as $\Pk$. In order to make this idea more clear we describe one particular
update rule, known as the \emph{heat-bath} update rule. This example
serves a dual purpose as it is a simple way to implement $\Pk$ and
we will make use of heat-bath updates in Sections~\ref{sec:gen2delta}
and~\ref{sec:tree} when applying our condition to specific spin
systems. A heat-bath move on a block $\Theta_{k}$ given a configuration
$x$ is performed by drawing a new configuration from the distribution
induced by $\pi$ and consistent with the assignment of spins on the
boundary of $\Theta_{k}$. The two properties of $\Pk$ hold for heat-bath
updates since (1) only the assignment of the spin to the sites in
$\Theta_{k}$ are changed and (2) the new configuration is drawn from
an appropriate distribution induced by $\pi$. If the spin system
corresponds to proper colourings of graphs then the distribution used
in a heat-bath move is the uniform distribution the set of configurations
that agree with $x$ off $\Theta_{k}$ and where no edge containing
a site in $\Theta_{k}$ is monochromatic. 

With these definitions in mind we are ready to formally define a systematic
scan Markov chain. 

\begin{defn}
We let $\scan$ be a systematic scan Markov chain with state space
$\Omega^{+}$ and transition matrix $\Pscan=\prod_{k=1}^{m}\Pk$. 
\end{defn}
The stationary distribution of $\scan$ is $\pi$ as discussed above,
and it is worth pointing out that the definition of $\scan$ holds
for \emph{any} order on the set of blocks. We will refer to one application
of $\Pscan$ (that is updating each block once) as one \emph{scan}
of $\scan$. One scan takes $\sum_{k}|\Theta_{k}|$ \emph{updates}
and it is generally straight forward to ensure, via the construction
of the set of blocks, that this sum is of order $O(n)$.

We will be concerned with analysing the mixing time of systematic
scan Markov chains, and consider the case when $\scan$ is ergodic.
Let $\M$ be any ergodic Markov chain with state space $\Omega^{+}$
and transition matrix $P$. By classical theory (see e.g. Aldous~\cite{aldous_walks})
$\M$ has a unique stationary distribution, which we will denote $\pi$.
The mixing time from an initial configuration $x\in\Omega^{+}$ is
the number of steps, that is applications of $P$, required for $\M$
to become sufficiently close to $\pi$. Formally the mixing time of
$\M$ from an initial configuration $x\in\Omega^{+}$ is defined,
as a function of the deviation $\epsilon$ from stationarity, by \[
\Mix_{x}(\M,\epsilon)=\min\{ t>0:\dtv(P^{t}(x,\cdot),\pi(\cdot))\leq\epsilon\}\]
where \[
\dtv(\theta_{1},\theta_{2})=\frac{1}{2}\sum_{i}|\theta_{1}(i)-\theta_{2}(i)|=\max_{A\subseteq\Omega^{+}}|\theta_{1}(A)-\theta_{2}(A)|\]
is the total variation distance between two distributions $\theta_{1}$
and $\theta_{2}$ on $\Omega^{+}$. The mixing time $\Mix(\M,\epsilon)$
of $\M$ is then obtained my maximising over all possible initial
configurations \[
\Mix(\M,\epsilon)=\max_{x\in\Omega^{+}}\Mix_{x}(\M,\epsilon).\]
 We say that $\M$ is \emph{rapidly mixing} if the mixing time of
$\M$ is polynomial in $n$ and $\log(\epsilon^{-1})$.

We will now formalise the notion of {}``the influence \emph{on} a
site'' in order to state our condition for rapid mixing of systematic
scan. For any pair of configurations $(x,y)$ let $\Psi_{k}(x,y)$
be a coupling of the distributions $\Pk(x)$ and $\Pk(y)$ which we
will refer to as {}``updating block $\Theta_{k}$''. Recall that
a coupling $\Psi_{k}(x,y)$ of $\Pk(x)$ and $\Pk(y)$ is a joint
distribution on $\Omega^{+}\times\Omega^{+}$ whose marginal distributions
are $\Pk(x)$ and $\Pk(y)$. That is \[
\forall\sigma\in\Omega^{+}\quad\Prob_{x^{\prime}\in\Pk(x)}(x^{\prime}=\sigma)=\sum_{\tau\in\Omega^{+}}\Prob_{(x^{\prime},y^{\prime})\in\Psi_{k}(x,y)}(x^{\prime}=\sigma,y^{\prime}=\tau)\]
and\[
\forall\tau\in\Omega^{+}\quad\Prob_{y^{\prime}\in\Pk(y)}(y^{\prime}=\sigma)=\sum_{\sigma\in\Omega^{+}}\Prob_{(x^{\prime},y^{\prime})\in\Psi_{k}(x,y)}(x^{\prime}=\sigma,y^{\prime}=\tau)\]
where we write $(x^{\prime},y^{\prime})\in\Psi_{k}(x,y)$ when the
pair of configurations $(x^{\prime},y^{\prime})$ is drawn from $\Psi_{k}(x,y)$.
Weitz in \cite{dror_combinatorial} states his conditions for general
metrics whereas we will use Hamming distance, which is also how the
corresponding condition is defined in Dyer et al. \cite{dobrushin_scan}.
This choice of metric allows us to define the influence of a site
$i$ on a site $j$ under a block $\Theta_{k}$, which we will denote
$\rho_{i,j}^{k}$, as the maximum probability that two coupled Markov
chains differ at the spin of site $j$ following an update of $\Theta_{k}$
starting from two configurations that only differ at the spin on site
$i$. That is \[
\rho_{i,j}^{k}=\max_{(x,y)\in S_{i}}\{\Prob_{(x^{\prime},y^{\prime})\in\Psi_{k}(x,y)}(x_{j}^{\prime}\neq y_{j}^{\prime})\}.\]
Then let $\alpha$ be the total (weighted) influence \emph{on} any
site in the graph site defined by \[
\alpha=\max_{k}\max_{j\in\Theta_{k}}\sum_{i}\frac{w_{i}}{w_{j}}\rho_{i,j}^{k}.\]
We point out that our definition of $\rho_{i,j}^{k}$ is not the standard
definition of $\rho$ used in the literature (see for example Simon
\cite{simon_lattice_gases} or Dyer et al. \cite{dobrushin_scan})
since the coupling $\Psi_{k}(x,y)$ is explicitly included. In the
block setting it is, however, necessary to include the coupling directly
in the definition of $\rho$ as we will discuss in Section \ref{sec:comparison-influence}.
In Section \ref{sec:comparison-influence} we also show that the condition
$\alpha<1$ is a generalisation of the corresponding condition in
Dyer et al. \cite{dobrushin_scan} in the sense that if each block
contains exactly one site and the coupling minimises the Hamming distance
then the conditions coincide. Our main theorem, which is proved in
Section \ref{sec:mixing}, states that if the influence on a site
is sufficiently small then the systematic scan Markov chain $\scan$
mixes in $O(\log n)$ scans. 

\begin{thm}
\label{thm:main_d} Suppose $\alpha<1$. Then \[
\Mix(\Mscan,\epsilon)\leq\frac{\log(n\epsilon^{-1})}{1-\alpha}.\]

\end{thm}
As previously stated we will apply Theorem \ref{thm:main_d} to two
spin systems corresponding to proper $q$-colourings of graphs in
order to improve the parameters for which systematic scan mixes. In
both applications we restrict the state space of the Markov chains
to the set of proper colourings, $\Omega$, of the underlying graph.
Firstly we allow the underlying graph to be any finite graph with
maximum vertex-degree $\Delta$. Previously, the least number of colours
for which systematic scan was known to mix in $O(\log n)$ scans was
$q>2\Delta$ and when $q=2\Delta$ the best known bound on the mixing
time was $O(n^{2}\log n)$ scans due to Dyer et al. \cite{dobrushin_scan}.
For completeness we pause to mention that the least number of colours
required for rapid mixing of a random update Markov chain is $q>11/6\Delta$
due to Vigoda \cite{vigoda}. In Section \ref{sec:gen2delta} we consider
the following Markov chain, \emph{edge scan} denoted $\escan$, updating
the endpoints of an edge during each update. Let $\Theta=\{\Theta_{k}\}_{k=1,\ldots,m}$
be a set of edges in $G$ such that $\Theta$ covers $V$. Using the
above notation, $\Pk$ is the transition matrix for performing a heat-bath
move on the endpoints of the edge $\Theta_{k}$ and the transition
matrix of $\escan$ is $\Pi_{k=1}^{m}\Pk$. We prove the following
theorem, which improves the mixing time of systematic scan by a factor
of $n^{2}$ for proper colourings of general graphs when $q=2\Delta$
and matches the existing bound when $q>2\Delta$. 

\begin{thm}
\label{thm:gen2delta} Let $G$ be a graph with maximum vertex-degree
$\Delta$. If $q\geq2\Delta$ then \[
\Mix(\escan,\epsilon)\leq\Delta^{2}\log(n\epsilon^{-1}).\]

\end{thm}
Next, in Section \ref{sec:tree}, we restrict the class of graphs
to trees. It is known that single site systematic scan mixes in $O(\log n)$
scans when $q>\Delta+2\sqrt{\Delta-1}$ and in $O(n^{2}\log n)$ scans
when $q=\Delta+2\sqrt{\Delta-1}$ is an integer; see e.g. Hayes \cite{hayes_dobrushin}
or Dyer, Goldberg and Jerrum \cite{matrix_norms}. More generally
it is known that systematic scan for proper colourings of bipartite
graphs mixes in $O(\log n)$ scans when $q\geq1.76\Delta$ as $\Delta\to\infty$
due to Bordewich et al. \cite{metric}. Again, for completeness, we
mention that the mixing time of a random update Markov chain for proper
colourings on a tree mixes in $O(n\log n)$ updates when $q\geq\Delta+2$,
a result due to Martinelli, Sinclair and Weitz \cite{tree}, improving
a similar result by Kenyon, Mossel and Peres~\cite{kenyon_tree}.
We will use a block approach to improve the number of colours required
for mixing of systematic scan on trees. We construct the following
set of blocks where the height $h$ of the blocks is defined in Table
\ref{tbl:block}. Let a block $\Theta_{k}$ contain a site $r$ along
with all sites below $r$ in the tree that are at most $h-1$ edges
away from $r$. The set of blocks $\Theta$ covers the sites of the
tree and we construct $\Theta$ such that no block has height less
than $h$. $\Pk$ is the transition matrix for performing a heat-bath
move on block $\Theta_{k}$ and the transition matrix of the Markov
chain $\treeblock$ is $\Pi_{k=1}^{m}\Pk$ where $m$ is the number
of blocks. We prove the following theorem. %
\begin{table}

\caption{Optimising the number of colours using blocks}

\begin{centering}\label{tbl:block} \begin{tabular}{c|c|c|c||c}
$\Delta$ &
$h$ &
$\xi$ &
$f(\Delta)$ &
$\lceil\Delta+2\sqrt{\Delta-1}\rceil$\tabularnewline
\hline 
3&
15&
$\frac{4}{7}$&
5&
6 \tabularnewline
4&
3&
$\frac{5}{11}$&
7&
8 \tabularnewline
5&
12&
$\frac{5}{11}$&
8&
9 \tabularnewline
6&
3&
$\frac{1}{2}$&
10&
11 \tabularnewline
7&
7&
$\frac{10}{23}$&
11&
12 \tabularnewline
8&
13&
$\frac{1}{3}$&
12&
14 \tabularnewline
9&
85&
$\frac{5}{19}$&
13&
15 \tabularnewline
10&
5&
$\frac{5}{19}$&
15&
16 \tabularnewline
\end{tabular}\par\end{centering}
\end{table}

\begin{thm}
\label{thm:tree_block} Let $G$ be a tree with maximum vertex-degree
$\Delta$. If $q\geq f(\Delta)$ where $f(\Delta)$ is specified in
Table $\ref{tbl:block}$ for small $\Delta$ then \[
\Mix(\treeblock,\epsilon)=O(\log(n\epsilon^{-1})).\]

\end{thm}
We conclude the paper with a discussion, in Section \ref{sec:comparison-influence},
of the influence parameter $\alpha$ and how it relates to the corresponding
parameters for the ``influence on a site'' in Weitz \cite{dror_combinatorial}
and Dyer et al. \cite{dobrushin_scan}. In particular we will show
that the condition in Weitz \cite{dror_combinatorial} does not imply
mixing of systematic scan and that the condition in Dyer et al. \cite{dobrushin_scan}
is a special case of our condition from Theorem \ref{thm:main_d}.

\section{Bounding the Mixing Time of Systematic Scan\label{sec:mixing} }

This section will contain the proof of Theorem \ref{thm:main_d}.
The proof follows the structure of the proof from the single-site
setting in Dyer et al. \cite{dobrushin_scan}, which follows F\"ollmer's
\cite{follmer} account of Dobrushin's proof presented in Simon's
book \cite{simon_lattice_gases}.

We will make use the following definitions. For any function $f:\Omega^{+}\rightarrow\R_{\geq0}$
let $\delta_{i}(f)=\max_{(x,y)\in S_{i}}|f(x)-f(y)|$ and $\Delta(f)=\sum_{i\in V}w_{i}\delta_{i}(f)$.
Also for any transition matrix $P$ define $(Pf)$ as the function
from $\Omega^{+}$ to $\R_{\geq0}$ given by $(Pf)(x)=\sum_{x^{\prime}}P(x,x^{\prime})f(x^{\prime})$.
Finally let $\mathbf{1}_{i\not\in\Theta_{k}}$ be the function given
by\[
\mathbf{1}_{i\not\in\Theta_{k}}=\begin{cases}
1 & \mbox{if }i\not\in\Theta_{k}\\
0 & \mbox{otherwise.}\end{cases}\]

We can think of $\delta_{i}(f)$ as the deviation from constancy of
$f$ at site $i$ and $\Delta(f)$ as the aggregated deviation from
constancy of $f$. Now, $Pf$ is a function where $(Pf)(x)$ gives
the expected value of $f$ after making a transition starting from
$x$. Intuitively, if $t$ transitions are sufficient for mixing then
$P^{t}f$ is a very smooth function. An application of $\Pk$ fixes
the non-constancy of $f$ at the sites within $\Theta_{k}$ although
possibly at the cost of increasing the non-constancy at sites on the
boundary of $\Theta_{k}$. Our aim is then to show that one application
of $\Pscan$ will on aggregate make $f$ smoother i.e., decrease $\Delta(f).$We
will establish the following lemma, which corresponds to Corollary
12 in Dyer et al. \cite{dobrushin_scan}, from which Section 3.3 of
\cite{dobrushin_scan} implies Theorem \ref{thm:main_d}.

\begin{lem}
\label{lemma:scan-sweep}If $\alpha<1$ then  \[
\Delta(\Pscan f)\leq\alpha\Delta(f).\]

\end{lem}
We begin by bounding the effect on $f$ from one application of $\Pk$.
The following lemma is a block-move generalisation of Proposition
V.1.7 from Simon \cite{simon_lattice_gases} and Lemma 10 from Dyer
et al. \cite{dobrushin_scan}. 

\begin{lem}
\label{lem:simon} $\delta_{i}(\Pk f)\leq\mathbf{1}_{i\not\in\Theta_{k}}\delta_{i}(f)+\sum_{j\in\Theta_{k}}\rho_{i,j}^{k}\delta_{j}(f)$ 
\end{lem}
\begin{proof}
Take $\Es{(x^{\prime},y^{\prime})\in\Psi_{k}(x,y))}{f(x^{\prime})}$
to be the the expected value of $f(x^{\prime})$ when a pair of configurations
$(x^{\prime},y^{\prime})$ are drawn from $\Psi_{k}(x,y)$. Since
$\Psi_{k}(x,y)$ is a coupling of the distributions $\Pk(x)$ and
$\Pk(y)$, the distribution $\Pk(x)$ and the first component of $\Psi_{k}(x,y)$
are the same and hence \begin{equation}
\Es{(x^{\prime},y^{\prime})\in\Psi_{k}(x,y)}{f(x^{\prime})}=\Es{x^{\prime}\in\Pk(x)}{f(x^{\prime})}\label{eq:expx}\end{equation}
and the same fact holds for the distribution $\Pk(y)$ so \begin{equation}
\Es{(x^{\prime},y^{\prime})\in\Psi_{k}(x,y)}{f(y^{\prime})}=\Es{y^{\prime}\in\Pk(y)}{f(y^{\prime})}.\label{eq:expy}\end{equation}
 Using \eqref{eq:expx}, \eqref{eq:expy} and linearity of expectation
we have \begin{align*}
\delta_{i}(\Pk f) & =\max_{(x,y)\in S_{i}}\left|(\Pk f)(x)-(\Pk f)(y)\right|\\
 & =\max_{(x,y)\in S_{i}}\left|\sum_{x^{\prime}}\Pk(x,x^{\prime})f(x^{\prime})-\sum_{y^{\prime}}\Pk(y,y^{\prime})f(y^{\prime})\right|\\
 & =\max_{(x,y)\in S_{i}}\left|\Es{x^{\prime}\in\Pk(x)}{f(x^{\prime})}-\Es{y^{\prime}\in\Pk(y)}{f(y^{\prime})}\right|\\
 & =\max_{(x,y)\in S_{i}}\left|\Es{(x^{\prime},y^{\prime})\in\Psi_{k}(x,y))}{f(x^{\prime})}-\Es{(x^{\prime},y^{\prime})\in\Psi_{k}(x,y)}{f(y^{\prime})}\right|\\
 & =\max_{(x,y)\in S_{i}}\left|\Es{(x^{\prime},y^{\prime})\in\Psi_{k}(x,y)}{f(x^{\prime})-f(y^{\prime})}\right|\\
 & \leq\max_{(x,y)\in S_{i}}\Es{(x^{\prime},y^{\prime})\in\Psi_{k}(x,y)}{\left|f(x^{\prime})-f(y^{\prime})\right|}\\
 & \leq\max_{(x,y)\in S_{i}}\Es{(x^{\prime},y^{\prime})\in\Psi_{k}(x,y)}{\sum_{j\in V}\left|f(x_{1}^{\prime}\dots x_{j}^{\prime}y_{j+1}^{\prime}\dots y_{n}^{\prime})-f(x_{1}^{\prime}\dots x_{j-1}^{\prime}y_{j}^{\prime}\dots y_{n}^{\prime})\right|}\\
 & =\max_{(x,y)\in S_{i}}\sum_{j\in V}\Es{(x^{\prime},y^{\prime})\in\Psi_{k}(x,y)}{\left|f(x_{1}^{\prime}\dots x_{j}^{\prime}y_{j+1}^{\prime}\dots y_{n}^{\prime})-f(x_{1}^{\prime}\dots x_{j-1}^{\prime}y_{j}^{\prime}\dots y_{n}^{\prime})\right|}.\end{align*}
 Notice that $x=x^{\prime}$ off $\Theta_{k}$ and $y=y^{\prime}$
off $\Theta_{k}$.

We need to bound the expectation $\Es{(x^{\prime},y^{\prime})\in\Psi_{k}(x,y)}{\left|f(x_{1}^{\prime}\dots x_{j}^{\prime}y_{j+1}^{\prime}\dots y_{n}^{\prime})-f(x_{1}^{\prime}\dots x_{j-1}^{\prime}y_{j}^{\prime}\dots y_{n}^{\prime})\right|}$
for each site $j\in V$. There are three cases. 
\begin{itemize}
\item $j\in\Theta_{k}$. By definition of $\rho_{i,j}^{k}$ the coupling
will yield $x_{j}^{\prime}\neq y_{j}^{\prime}$ with probability at
most $\rho_{i,j}^{k}$ and so \begin{eqnarray*}
 &  & \Es{(x^{\prime},y^{\prime})\in\Psi_{k}(x,y)}{\left|f(x_{1}^{\prime}\dots x_{j}^{\prime}y_{j+1}^{\prime}\dots y_{n}^{\prime})-f(x_{1}^{\prime}\dots x_{j-1}^{\prime}y_{j}^{\prime}\dots y_{n}^{\prime})\right|}\\
 &  & \quad\leq\rho_{i,j}^{k}\max_{(\sigma,\tau)\in S_{j}}\{|f(\sigma)-f(\tau)|\}=\rho_{i,j}^{k}\delta_{j}(f).\end{eqnarray*}
 
\item $j\not\in\Theta_{k}$ and $j=i$. Since $j\not\in\Theta_{k}$ we have
$x_{j}=x_{j}^{\prime}$ and $y_{j}=y_{j}^{\prime}$ so \[
\Es{(x^{\prime},y^{\prime})\in\Psi_{k}(x,y)}{\left|f(x_{1}^{\prime}\dots x_{j}^{\prime}y_{j+1}^{\prime}\dots y_{n}^{\prime})-f(x_{1}^{\prime}\dots x_{j-1}^{\prime}y_{j}^{\prime}\dots y_{n}^{\prime})\right|}\leq\delta_{j}(f)=\delta_{i}(f).\]
 
\item $j\not\in\Theta_{k}$ and $i\neq j$. In this case we have $x_{j}=x_{j}^{\prime}$
and $y_{j}=y_{j}^{\prime}$ which implies $x_{j}^{\prime}=y_{j}^{\prime}$
so \[
\Es{(x^{\prime},y^{\prime})\in\Psi_{k}(x,y)}{\left|f(x_{1}^{\prime}\dots x_{j}^{\prime}y_{j+1}^{\prime}\dots y_{n}^{\prime})-f(x_{1}^{\prime}\dots x_{j-1}^{\prime}y_{j}^{\prime}\dots y_{n}^{\prime})\right|}=0.\]

\end{itemize}
Adding it up we get the statement of the lemma. 
\end{proof}
We will use Lemma \ref{lem:simon} in conjunction with an inductive
proof similar to (V.1.16) in Simon \cite{simon_lattice_gases} in
order to establish the following lemma. It is important to note at
this point that the result in Simon is presented for single site heat-bath
updates, whereas the following lemma applies to any block dynamics
(satisfying the stated assumptions) and weighted sites. This lemma
is also a block generalisation of Lemma 11 in Dyer et al. \cite{dobrushin_scan}.

\begin{lem}
\label{lemma:k_scan} Let $\Gamma(k)=\bigcup_{l=1}^{k}\Theta_{l}$
then for any $k\in\{1,\dots,m\}$, if $\alpha<1$ then \[
\Delta(\P1\cdots\Pk f)\leq\alpha\sum_{i\in\Gamma(k)}w_{i}\delta_{i}(f)+\sum_{i\in V\setminus\Gamma(k)}w_{i}\delta_{i}(f).\]

\end{lem}
\begin{proof}
Induction on $k$. Taking $k=0$ as the base case, we get the definition
of $\Delta$.

Assume the statement holds for $k-1$. \[
\begin{split}\Delta(\P1\cdots\Pk f) & \leq\alpha\sum_{i\in\Gamma(k-1)}w_{i}\delta_{i}(\Pk f)+\sum_{i\in V\setminus\Gamma(k-1)}w_{i}\delta_{i}(\Pk f)\\
 & \leq\alpha\sum_{i\in\Gamma(k-1)}\mathbf{1}_{i\not\in\Theta_{k}}w_{i}\delta_{i}(f)+\alpha\sum_{i\in\Gamma(k-1)}\sum_{j\in\Theta_{k}}w_{i}\rho_{i,j}^{k}\delta_{j}(f)\\
 & \quad+\sum_{i\in V\setminus\Gamma(k-1)}\mathbf{1}_{i\not\in\Theta_{k}}w_{i}\delta_{i}(f)+\sum_{i\in V\setminus\Gamma(k-1)}\sum_{j\in\Theta_{k}}w_{i}\rho_{i,j}^{k}\delta_{j}(f)\end{split}
\]
 by Lemma \ref{lem:simon}.

Simplifying and using $\alpha<1$ \[
\begin{split}\Delta(\P1\cdots\Pk f) & \leq\alpha\sum_{i\in\Gamma(k-1)\setminus\Theta_{k}}w_{i}\delta_{i}(f)+\sum_{i\in\Gamma(k-1)}\sum_{j\in\Theta_{k}}w_{i}\rho_{i,j}^{k}\delta_{j}(f)\\
 & \quad+\sum_{i\in V\setminus\Gamma(k)}w_{i}\delta_{i}(f)+\sum_{i\in V\setminus\Gamma(k-1)}\sum_{j\in\Theta_{k}}w_{i}\rho_{i,j}^{k}\delta_{j}(f)\\
 & =\alpha\sum_{i\in\Gamma(k-1)\setminus\Theta_{k}}w_{i}\delta_{i}(f)+\sum_{i\in V\setminus\Gamma(k)}w_{i}\delta_{i}(f)\\
 & \quad+\sum_{j\in\Theta_{k}}\delta_{j}(f)\left(\sum_{i\in\Gamma(k-1)}w_{i}\rho_{i,j}^{k}+\sum_{i\in V\setminus\Gamma(k-1)}w_{i}\rho_{i,j}^{k}\right)\\
 & =\alpha\sum_{i\in\Gamma(k-1)\setminus\Theta_{k}}w_{i}\delta_{i}(f)+\sum_{i\in V\setminus\Gamma(k)}w_{i}\delta_{i}(f)+\sum_{j\in\Theta_{k}}\delta_{j}(f)\sum_{i\in V}w_{i}\rho_{i,j}^{k}\\
 & \leq\alpha\sum_{i\in\Gamma(k-1)\setminus\Theta_{k}}w_{i}\delta_{i}(f)+\sum_{i\in V\setminus\Gamma(k)}w_{i}\delta_{i}(f)+\sum_{j\in\Theta_{k}}\delta_{j}(f)\max_{l}\sum_{i\in V}w_{i}\rho_{i,j}^{l}\\
 & \leq\alpha\sum_{i\in\Gamma(k-1)\setminus\Theta_{k}}w_{i}\delta_{i}(f)+\sum_{i\in V\setminus\Gamma(k)}w_{i}\delta_{i}(f)+\alpha\sum_{j\in\Theta_{k}}w_{j}\delta_{j}(f)\\
 & =\alpha\sum_{i\in\Gamma(k)}w_{i}\delta_{i}(f)+\sum_{i\in V\setminus\Gamma(k)}w_{i}\delta_{i}(f)\end{split}
\]
 by definition of $\alpha$.
\end{proof}
Lemma \ref{lemma:scan-sweep} is now a simple consequence of Lemma
\ref{lemma:k_scan} since \[
\Delta(\Pscan f)=\Delta(\P1\cdots\Pm f)\leq\alpha\sum_{i\in V}w_{i}\delta_{i}(f)=\alpha\Delta(f)\]
and Theorem \ref{thm:main_d} follows as discussed above.

\section{Application: Edge Scan on an Arbitrary Graph\label{sec:gen2delta}}

In this section we prove Theorem \ref{thm:gen2delta}. That is, we
present a general version of a systematic scan on edges and use Theorem
\ref{thm:main_d} to prove that it mixes in $O(\log n)$ scans when
$q\geq2\Delta$. We use $w_{i}=1$ for all $i\in V$ and so omit all
weights throughout this section. Recall that $\escan$ is the systematic
scan Markov chain with transition matrix $\Pi_{k=1}^{m}\Pk$ where
$\Theta=\{\Theta_{k}\}_{k=1,\ldots,m}$ is an ordered set of edges
in $G$ that covers $V$ and $\Pk$ is the transition matrix for performing
a heat-bath move on the endpoints of the edge $\Theta_{k}$.

We need to construct a coupling $\Psi_{k}(x,y)$ of the distributions
$\Pk(x)$ and $\Pk(y)$ for each pair of configurations $(x,y)\in S_{i}$
that differ only at the colour assigned to site $i$. Assume without
loss of generality that $x_{i}=1$ and $y_{i}=2$ and also let $j$
and $j^{\prime}$ be the endpoints of the edge $\Theta_{k}$. Recall
that, since the dynamics uses heat-bath updates, $\Pk(x)$ is the
uniform distribution on configurations that agree with $x$ off $\Theta_{k}$
and where no edge containing $j$ or $j^{\prime}$ is monochromatic.
For ease of notation we let $D_{1}=\Pk(x)$ and $D_{2}=\Pk(y)$. We
go on to make the following definitions for $l\in\{1,2\}$ and $s\in\Theta_{k}$.
$D_{l}(s)$ is the distribution of the colour assigned to site $s$
induced by $D_{l}$, and $[D{}_{l}\mid s=c]$ is the uniform distribution
on the set of colourings of the sites in $\Theta_{k}$ where site
$s$ is assigned colour $c$. We also let $d_{l}$ denote the number
of configurations with positive measure in $D_{l}$ and $d_{l,s=c}$
be the number of configurations that assign colour $c$ to site $s$
and have positive measure in $D_{l}$.

\begin{defn}
We will say  that the choice $c_{1}c_{2}$ is ``valid'' for $D_{l}$
if there is a configuration with positive measure in $D_{l}$ in which
site $j$ is coloured $c_{1}$ and site $j^{\prime}$ is coloured
$c_{2}$. Similarly a colour $c$ is {}``valid'' on a site $s$
in $D_{l}$ if there exists a valid choice for $D_{l}$ where site
$s$ is coloured $c$. 
\end{defn}

\subsection{Overview of the Coupling}

We begin the construction of the coupling $\Psi_{k}(x,y)$ by giving
an overview of the cases we will need to consider and show that they
are mutually exclusive and exhaustive of all configurations. It is
important to note that, by definition of $\rho$, the coupling we
define may depend on the initial configurations $x$ and $y$ in the
sense that if two pairs of configurations $(x_{1},y_{1})$ and $(x_{2},y_{2})$
can be distinguished then the couplings $\Psi_{k}(x_{1},y_{1})$ and
$\Psi_{k}(x_{2},y_{2})$ may be defined differently. 

First, if $i$ is not adjacent to any site in $\Theta_{k}$, that
is $i\not\in\partial\Theta_{k}$, then $\Psi_{k}(x,y)$ is the identity
coupling where the same colouring is assigned to each distribution.
Hence, for $i\not\in\partial\Theta_{k}$ and $j\in\Theta_{k}$ we
have\[
\rho_{i,j}^{k}=0.\]

Now suppose that $i$ is adjacent to at least one site in $\Theta_{k}$,
that is $i\in\partial\Theta_{k}$. We consider the following five
cases, which by construction are exhaustive of all possible configurations
and mutually exclusive. In the diagrams that relate to these cases
a dotted line between a site $j\in\Theta_{k}$ and a colour $1$,
say, denotes that no site adjacent to $j$ on the boundary of $\Theta_{k}$
(other than possibly $i$) is coloured $1$. A full line denotes that
some site adjacent to $j$ on the boundary of $\Theta_{k}$ (other
than possibly $i$) is coloured $1$. The full details of each case
of the coupling will be given in section \ref{subsec:fullcup} along
with bounds on $\rho_{i,j}^{k}$ and $\rho_{i,j^{\prime}}^{k}$ where
$j$ and $j^{\prime}$ are the sites included in $\Theta_{k}$.

\def\one{Exactly one site in $\Theta_k$ is adjacent to $i$. Let this site be labeled $j$ and let the other site in $\Theta_k$ be labeled $j^\prime$}
\def\three{Both sites in $\Theta_k$ are adjacent to $i$ and no other sites in $\partial \Theta_k$ are coloured 1 or 2. The labeling of the sites in $\Theta_k$ is arbitrary}
\def\four{Both sites in $\Theta_k$ are adjacent to $i$. One of the sites in $\Theta_k$ is adjacent to at least one site, other than $i$, coloured 1. Let this site be labeled $j^\prime$. The other site in $\Theta_k$ is labeled $j$ and it is not adjacent to any site, other than $i$, coloured 1 or 2}
\def\five{Both sites in $\Theta_k$ are adjacent to $i$. One of the sites in $\Theta_k$ is adjacent to at least one site, other than $i$, coloured 1 and no sites that are coloured 2. Let this site be labeled $j^\prime$. The other site in $\Theta_k$, labeled $j$, is adjacent to at least one site other than $i$ coloured 2 and no sites coloured 1}
\def\six{Both sites in $\Theta_k$ are adjacent to $i$ and at least one site, other than $i$ coloured 1. The labeling of the sites in $\Theta_k$ is arbitrary}

\begin{figure}

\caption{Case \ref{c1}. \one.}

\label{fig:c1} 

\begin{centering}\includegraphics{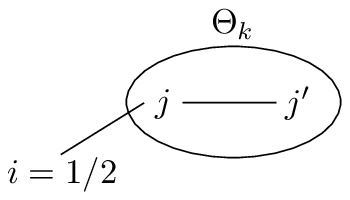} \par\end{centering}
\end{figure}
\begin{figure}

\caption{Case \ref{c3}. \three.}

\label{fig:c3} 

\begin{centering}\includegraphics{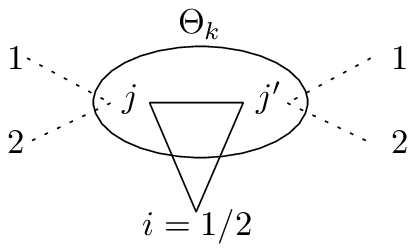} \par\end{centering}
\end{figure}
\begin{figure}

\caption{Case \ref{c4}. \four.}

\label{fig:c4} 

\begin{centering}\includegraphics{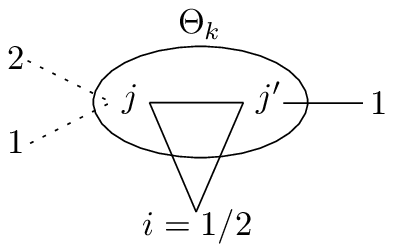} \par\end{centering}
\end{figure}
\begin{figure}

\caption{Case \ref{c5}. \five.}

\label{fig:c5} 

\begin{centering}\includegraphics{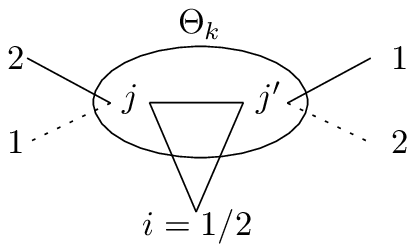} \par\end{centering}
\end{figure}
\begin{figure}

\caption{Case \ref{c6}. \six.}

\label{fig:c6} 

\begin{centering}\includegraphics{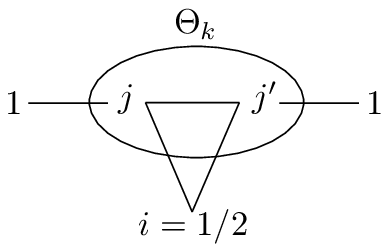} \par\end{centering}
\end{figure}

\begin{enumerate}
\item \label{c1}\one. This is shown in Figure \ref{fig:c1}. 
\item \label{c3}\three. This is shown in Figure \ref{fig:c3}. 
\item \label{c4}\four. This is shown in Figure \ref{fig:c4}. 
\item \label{c5}\five. This is shown in Figure \ref{fig:c5}. 
\item \label{c6}\six. This is shown in Figure \ref{fig:c6}.
\end{enumerate}

\subsection{Details of Coupling and Proof of Mixing\label{subsec:fullcup}}

We will now give the full details of each case of the coupling and
establish the required bounds on the influence of site $i$ on sites
$j$ and $j^{\prime}$. The following lemma is required to establish
the coupling for all the stated cases. 

\begin{figure}

\caption{Case \ref{c1}. \one.}

\label{fig:c1p} 

\begin{centering}\includegraphics{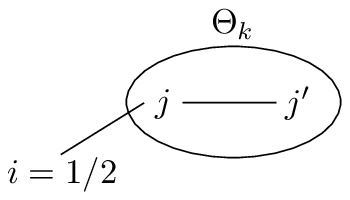} \par\end{centering}
\end{figure}

\begin{lem}
\label{lemma:d_same} Let $j$ and $j^{\prime}$ be the endpoints
of an edge $\Theta_{k}$ and suppose that $\{ i,j\}\in E(G)$. Then
for each pair of colours $c_{1},c_{2}\in C\setminus\{1,2\}$ the choice
$c_{1}c_{2}\textnormal{ is valid for }D_{1}$ if and only if $c_{1}c_{2}\textnormal{ is valid for }D_{2}.$
\end{lem}
\begin{proof}
We start with the \emph{if} direction. Suppose $c_{1}c_{2}$ is valid
in $D_{2}$ then no site adjacent to $j$ has colour $c_{1}$ in $D_{2}$
and since $c_{1}\neq1$ no site adjacent to $j$ has colour $c_{1}$
in $D_{1}$. Also no site adjacent to $j^{\prime}$ has colour $c_{2}$
in $D_{2}$ hence no site adjacent to $j^{\prime}$ has colour $c_{2}$
in $D_{1}$ since $c_{2}\neq1$. Since $c_{1}c_{2}$ is valid in $D_{2}$
$c_{1}\neq c_{2}$ and so $c_{1}c_{2}$ is valid in $D_{1}$.

The \emph{only if} direction is similar. Suppose $c_{1}c_{2}$ is
valid in $D_{1}$ then no site adjacent to $j$ has colour $c_{1}$
in $D_{1}$ and since $c_{1}\neq2$ no site adjacent to $j$ has colour
$c_{1}$ in $D_{2}$. Also no site adjacent to $j^{\prime}$ has colour
$c_{2}$ in $D_{1}$ hence no site adjacent to $j^{\prime}$ has colour
$c_{2}$ in $D_{2}$ again since $c_{2}\neq2$. Since $c_{1}c_{2}$
is valid in $D_{1}$ $c_{1}\neq c_{2}$ and so $c_{1}c_{2}$ is valid
in $D_{2}$. 
\end{proof}
\textbf{Details of case \ref{c1}.} (Repeated in Figure \ref{fig:c1p})
 We construct a coupling $\Psi_{k}(x,y)$ of the distributions $D_{1}$
and $D_{2}$ using the following two step process. Let $\psi_{j}$
be a coupling of $D_{1}(j)$ and $D_{2}(j)$ which greedily maximises
the probability of assigning the same colour to site $j$ in each
distribution. Then, for each pair of colours $(c,c^{\prime})$ drawn
from $\psi_{j}$, $\Psi_{k}(x,y)$ is a coupling, minimising Hamming
distance, of the conditional distributions $D_{1}\mid j=c$ and $D_{2}\mid j=c^{\prime}$. 

\begin{lem}
\label{lemma:gen2delta_notriag} Let $j$ and $j^{\prime}$ be the
endpoints of an edge $\Theta_{k}$. If $\{ i,j\}\in E(G)$ and $\{ i,j^{\prime}\}\not\in E(G)$
then\begin{align*}
\rho_{i,j}^{k}\leq\frac{1}{q-\Delta}\textnormal{ and }\rho_{i,j^{\prime}}^{k}\leq\frac{1}{(q-\Delta)^{2}}.\end{align*}

\end{lem}
\begin{proof}
Assume without loss of generality that $d_{1}\geq d_{2}$, i.e that
there are at least as many valid choices for $D_{1}$ as for $D_{2}$.
Since the only site in $\Theta_{k}$ that is adjacent to site $i$
is $j$, Lemma 13 of Goldberg, Martin and Paterson \cite{ssm} lets
us upper bound the probability of a discrepancy at site $j$ in a
pair of configurations drawn from the coupling $\Psi_{k}(x,y)$ by
assuming that $j^{\prime}$ is assigned the worst case colour. Now,
$1$ is not valid for $j$ in $D_{1}$ so Lemma \ref{lemma:d_same}
implies that only the choice $2$ for $j$ in $D_{1}$ would cause
site $j$ to be assigned a different colour in each configuration
drawn from the coupling. Now observe that site $j$ has at most $\Delta-1$
neighbours (excluding $j^{\prime}$) and each of them could invalidate
one colour choice for $j$ in both distributions. If $j^{\prime}$
is assigned a colour not already adjacent to $j$ then $j$ is adjacent
to at most $\Delta$ sites each assigned a different colour, leaving
at least $q-\Delta$ valid colours for $j$ in $D_{1}$ and so the
probability of assigning $2$ to $j$ in $D_{1}$ during step $1$
of the coupling is at most $\frac{1}{q-\Delta}$ since the coupling
is greedy. This establishes the bound on $\rho_{i,j}^{k}$ since\[
\rho_{i,j}^{k}=\max_{(x,y)\in S_{i}}\{\Prob_{(x^{\prime},y^{\prime})\in\Psi_{k}(x,y)}(x_{j}^{\prime}\neq y_{j}^{\prime})\}\leq\frac{1}{q-\Delta}.\]

Now from the definition of the coupling it follows easily that if
the same colour, $c,$ is assigned to site $j$ in each distribution
during the first step of the coupling then the colour assigned to
site $j^{\prime}$ in the second step will be the same in each distribution
since the conditional distributions $D_{1}\mid j=c$ and $D_{2}\mid j=c$
are the same. If different colours are assigned to $j$ in each distribution
then the second step of the coupling is simply the case of colouring
a single site adjacent to exactly one discrepancy. The argument from
above says that at most one colour assigned to $j^{\prime}$ in $D_{1}$
will cause a discrepancy at site $j^{\prime}$ in the coupling and
also that there are at least $q-\Delta$ valid choices for $j^{\prime}$
in $D_{1}$. Hence we have $\max_{(x,y)\in S_{i}}\{\Prob_{(x^{\prime},y^{\prime})\in\Psi_{k}(x,y)}(x_{j^{\prime}}^{\prime}\neq y_{j^{\prime}}^{\prime}\mid x_{j}^{\prime}=c,y_{j}^{\prime}=c^{\prime})\}\leq\frac{1}{q-\Delta}$
and so \begin{align*}
\rho_{i,j^{\prime}}^{k} & =\max_{(x,y)\in S_{i}}\{\Prob_{(x^{\prime},y^{\prime})\in\Psi_{k}(x,y)}(x_{j^{\prime}}^{\prime}\neq y_{j^{\prime}}^{\prime})\}\\
 & =\max_{(x,y)\in S_{i}}\left\{ \sum_{c,c^{\prime}}\Prob_{(x^{\prime},y^{\prime})\in\Psi_{k}(x,y)}(x_{j^{\prime}}^{\prime}\neq y_{j^{\prime}}^{\prime}\mid x_{j}^{\prime}=c,y_{j}^{\prime}=c^{\prime})\Prob_{(x^{\prime},y^{\prime})\in\Psi_{k}(x,y)}(x_{j}^{\prime}=c,y_{j}^{\prime}=c^{\prime})\right\} \\
 & \leq\frac{1}{q-\Delta}\max_{(x,y)\in S_{i}}\left\{ \sum_{c,c^{\prime}}\Prob_{(x^{\prime},y^{\prime})\in\Psi_{k}(x,y)}(x_{j}^{\prime}=c,y_{j}^{\prime}=c^{\prime})\right\} \\
 & \leq\frac{1}{(q-\Delta)^{2}}\end{align*}
 using the bound from $\rho_{i,j}^{k}$ which completes the proof.
\end{proof}
The following lemmas are required to define the coupling and bound
the influence of a site $i\in\partial\Theta_{k}$ on sites $j$ and
$j^{\prime}$ when $i$ is adjacent to both sites $j$ and $j^{\prime}$. 

\begin{lem}
\label{lemma:d} Let $j$ and $j^{\prime}$ be the endpoints of an
edge and suppose that $\{ i,j\}\in E(G)$ and $\{ i,j^{\prime}\}\in E(G)$.
If $1$ is valid for $j$ in $D_{2}$ and $2$ is valid for $j$ in
$D_{1}$ then the choice $2c_{2}\textnormal{ is valid in }D_{1}$
if and only if $1c_{2}\textnormal{ is valid in }D_{2}.$
\end{lem}
\begin{proof}
Suppose that $2c_{2}$ is valid in $D_{1}$ then $c_{2}\in C\setminus\{1,2\}$
since $i$ is adjacent to $j^{\prime}$ (and $x_{i}=1$). Since $1$
is valid for $j$ in $D_{2}$ it follows that $1c_{2}$ is valid in
$D_{2}$ since the only colour adjacent to $j^{\prime}$ in $D_{2}$
that is (possibly) not adjacent to $j^{\prime}$ in $D_{1}$ is $2$,
but $c_{2}\neq2$.

For the reverse direction suppose that $1c_{2}$ is valid in $D_{2}$.
Then $c_{2}\in C\setminus\{1,2\}$ since $i$ is adjacent to $j^{\prime}$.
Since $2$ is valid for $j$ in $D_{1}$ it follows that $2c_{2}$
is valid in $D_{1}$ since the only colour adjacent to $j^{\prime}$
in $D_{1}$ that is (possibly) not adjacent to $j^{\prime}$ in $D_{2}$
is $1$, but $c_{2}\neq1$.
\end{proof}
\begin{lem}
\label{lemma:dprime} Let $j$ and $j^{\prime}$ be the endpoints
of an edge $\Theta_{k}$ and suppose that $\{ i,j\}\in E(G)$ and
$\{ i,j^{\prime}\}\in E(G)$. If $1$ is valid for $j^{\prime}$ in
$D_{2}$ and $2$ is valid for $j^{\prime}$ in $D_{1}$ then the
choice $c_{1}2\textnormal{ is valid in }D_{1}$ if and only if $c_{1}1\textnormal{ is valid in }D_{2}.$
\end{lem}
\begin{proof}
Suppose that $c_{1}2$ is valid in $D_{1}$ then $c_{1}\in C\setminus\{1,2\}$
since $i$ is adjacent to $j^{\prime}$. Since $1$ is valid for $j^{\prime}$
in $D_{2}$ $c_{1}1$ is valid in $D_{2}$ since the only colour adjacent
to $j$ in $D_{2}$ that is (possibly) not adjacent to $j$ in $D_{1}$
is $2$, but $c_{1}\neq2$.

Also, suppose that $c_{1}1$ is valid in $D_{2}$ then $c_{1}\in C\setminus\{1,2\}$
since $i$ is adjacent to $j^{\prime}$. Since $2$ is valid for $j^{\prime}$
in $D_{1}$ $c_{1}2$ is valid in $D_{1}$ since the only colour adjacent
to $j$ in $D_{1}$ that is (possibly) not adjacent to $j$ in $D_{2}$
is $1$, but $c_{1}\neq1$. 
\end{proof}
\begin{lem}
\label{lemma:triangle_d} Let $j$ and $j^{\prime}$ be the endpoints
of an edge $\Theta_{k}$ and suppose that $\{ i,j\}\in E(G)$ and
$\{ i,j^{\prime}\}\in E(G)$.

\begin{enumerate}[{(i)}]
\item Suppose that $1$ is valid for $j$ in $D_{2}$. For all $c\in C$
where $c$ is valid for $j$ in $D_{2}$, if $1$ is valid for $j^{\prime}$
in $D_{2}$ then \[
d_{2,j=1}\leq d_{2,j=c}\leq d_{2,j=1}+1\]
 else \[
d_{2,j=1}-1\leq d_{2,j=c}\leq d_{2,j=1}.\]

\item Suppose that $2$ is valid for $j$ in $D_{1}$. For all $c\in C$
where $c$ is valid for $j$ in $D_{1}$, if $2$ is valid for $j^{\prime}$
in $D_{1}$ then \[
d_{1,j=2}\leq d_{1,j=c}\leq d_{1,j=2}+1\]
 else \[
d_{1,j=2}-1\leq d_{1,j=c}\leq d_{1,j=2}.\]
 \end{enumerate}
\end{lem}
\begin{proof}
Part (i). Consider some valid colour $c$ other than $1$ for $j$
in $D_{2}$. For each valid choice $1c_{2}$ for $D_{2}$ the choice
$cc_{2}$ is also valid for $D_{2}$ except when $c=c_{2}$. If $1$
is valid for $j^{\prime}$ in $D_{2}$ then the choice $c1$ is also
valid for $D_{2}$.

Now consider some invalid choice $1c_{2}$ for $D_{2}$ where $c_{2}\neq1$.
Since $1c_{2}$ is not valid for $D_{2}$ it follows that $c_{2}$
is not valid for $j^{\prime}$ in $D_{2}$ and hence no more choices
can be valid for $D_{2}$, which guarantees the upper bounds.

Part (ii) is similar. Consider some valid colour $c$ other than $2$
for $j$ in $D_{1}$. For each valid choice $2c_{2}$ for $D_{1}$
the choice $cc_{2}$ is also valid for $D_{1}$ except when $c=c_{2}$.
If $2$ is valid for $j^{\prime}$ in $D_{1}$ then the choice $c2$
is also valid for $D_{1}$.

Finally consider some invalid choice $2c_{2}$ for $D_{1}$ where
$c_{2}\neq2$. Since $2c_{2}$ is not valid for $D_{1}$ it follows
that $c_{2}$ is not valid for $j^{\prime}$ in $D_{1}$ and hence
no more choices can be valid for $D_{1}$, which guarantees the upper
bounds.
\end{proof}
We are now ready to define the coupling for the remaining cases.

\begin{figure}

\caption{Case \ref{c3}. \three.}

\label{fig:c3p} 

\begin{centering}\includegraphics{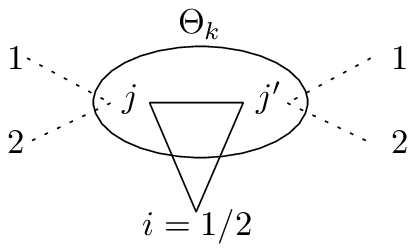} \par\end{centering}
\end{figure}

\textbf{Details of case \ref{c3}.} (Repeated in Figure \ref{fig:c3p})
 We construct the $\Psi_{k}(x,y)$ of the distributions $D_{1}$ and
$D_{2}$ as follows. For each valid choice of the form $c_{1}c_{2}$
for $D_{1}$ where $c_{1}\neq2$ and $c_{2}\neq2$ Lemma \ref{lemma:d_same}
guarantees that $c_{1}c_{2}$ is valid for $D_{2}$ so we let\[
\Prob_{(x^{\prime},y^{\prime})\in\Psi_{k}(x,y)}(x^{\prime}=y^{\prime}=c_{1}c_{2})=\frac{1}{d_{1}}.\]
For each valid choice of the form $2c_{2}$ in $D_{1}$ the choice
$1c_{2}$ is valid in $D_{2}$ by Lemma \ref{lemma:d} so we let\begin{equation}
\Prob_{(x^{\prime},y^{\prime})\in\Psi_{k}(x,y)}(x^{\prime}=2c_{2},y^{\prime}=1c_{2})=\frac{1}{d_{1}}.\label{eq:c3dis1}\end{equation}
Lemma \ref{lemma:d} also guarantees that there are no remaining valid
choices for $D_{2}$ of the form $1c_{2}$. Finally for each valid
choice $c_{1}2$ for $D_{1}$ the choice $c_{1}1$ is valid in $D_{2}$
by Lemma \ref{lemma:dprime} so let \begin{equation}
\Prob_{(x^{\prime},y^{\prime})\in\Psi_{k}(x,y)}(x^{\prime}=c_{1}2,y^{\prime}=c_{1}1)=\frac{1}{d_{1}}\label{eq:c3dis2}\end{equation}
which completes the coupling since $d_{1}=d_{2}$ and all the probability
in both $D_{1}$ and $D_{2}$ has hence been used.

\begin{figure}

\caption{Case \ref{c4}. \four.}

\label{fig:c4p} 

\begin{centering}\includegraphics{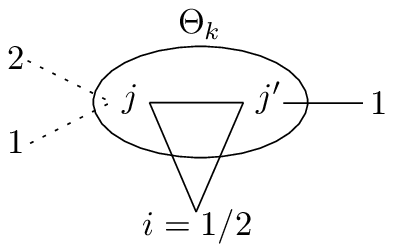} \par\end{centering}
\end{figure}

\begin{lem}
\label{lemma:tr1} Let $j$ and $j^{\prime}$ be the endpoints of
an edge $\Theta_{k}$ and suppose that $\{ i,j\}\in E(G)$ and $\{ i,j^{\prime}\}\in E(G)$.
If $2$ is valid for both $j$ and $j^{\prime}$ in $D_{1}$ and $1$
is valid for both $j$ and $j^{\prime}$ in $D_{2}$ then \[
\rho_{i,j}^{k}\leq\frac{1}{q-\Delta+1}\textnormal{ and }\rho_{i,j^{\prime}}^{k}\leq\frac{1}{q-\Delta}.\]

\end{lem}
\begin{proof}
This is case \ref{c3} of the coupling. Note from Lemma \ref{lemma:d}
that $d_{1,j=2}=d_{2,j=1}$ so for ease of reference let $d=d_{1,j=2}=d_{2,j=1}$
and let $d^{\prime}=d_{1,j^{\prime}=2}=d_{2,j^{\prime}=1}$ by Lemma
\ref{lemma:dprime}. Also let $s=\sum_{c}d_{2,j=c}-d-d^{\prime}$
which is the number of valid choices for $D_{2}$ other than choices
of the form $1c_{2}$ and $c_{1}1$. Note that the number of valid
choices for $D_{1}$ is $d_{1}=s+d+d^{\prime}$.

As there are no restrictions on colours assigned to the sites in $\partial\Theta_{k}\setminus\{ i\}$
each of the neighbours of $j$ could be assigned a different colour,
and the same is true for the neighbours of $j^{\prime}$. Hence we
get the following lower-bounds on $d$ and $d^{\prime}$: \[
q-\Delta\leq d\textnormal{ and }q-\Delta\leq d^{\prime}.\]
 To lower bound bound $s$ observe that $s=\sum_{c}d_{2,j=c}-d-d^{\prime}=\sum_{c\neq1}d_{2,j=c}-d^{\prime}$.
Let $J\subseteq C\setminus\{1\}$ be the set of colours, excluding
$1$, that are valid for $j$ in $D_{2}$. By definition of $d^{\prime}$,
at least $d^{\prime}$ colours other than $1$ must be valid for site
$j$ in $D_{2}$ so the size of $J$ is at least $d^{\prime}$. Since
$1$ is valid for $j^{\prime}$ in $D_{2}$ we use the lower bound
on $d_{2,j=c}$ from Lemma \ref{lemma:triangle_d} $(i)$ and hence\begin{align*}
s & =\sum_{c\in J}d_{2,j=c}-d^{\prime}\\
 & \geq d^{\prime}\min_{c\in J}\{ d_{2,j=c}\}-d^{\prime}\\
 & \geq d^{\prime}d-d^{\prime}.\end{align*}
From the coupling, $j$ will be assigned a different colour in each
distribution whenever a choice of the form $2c_{2}$ is made for $D_{1}$.
From \eqref{eq:c3dis1} this happens with probability $\frac{d}{d_{1}}=\frac{d}{d+d^{\prime}+s}$
since $d$ is the number of valid choices for $D_{1}$ of the form
$2c_{2}$. Similarly from \eqref{eq:c3dis2}, $j^{\prime}$ will become
a discrepancy in the coupling whenever a choice of the form $c_{1}2$
is made for $D_{1}$, which happens with probability $\frac{d^{\prime}}{d+d^{\prime}+s}$.
Hence \[
\rho_{i,j}^{k}\leq\frac{d}{d+d^{\prime}+s}\textnormal{ and }\rho_{i,j^{\prime}}^{k}\leq\frac{d^{\prime}}{d+d^{\prime}+s}.\]
 Starting with $\rho_{i,j}^{k}$ \[
\rho_{i,j}^{k}\leq\frac{d}{d+d^{\prime}+s}\leq\frac{d}{d+dd^{\prime}}\leq\frac{1}{d^{\prime}+1}\leq\frac{1}{q-\Delta+1}\]
 using the lower bounds of $s$ and $d^{\prime}$. Similarly using
the lower bounds of $s$ and $d$ \[
\rho_{i,j^{\prime}}^{k}\leq\frac{d^{\prime}}{d+d^{\prime}+s}\leq\frac{d^{\prime}}{d+dd^{\prime}}\leq\frac{1}{d}\leq\frac{1}{q-\Delta}\]
 which implies the statement of the lemma.
\end{proof}
\textbf{Details of case \ref{c4}.} (Repeated in Figure \ref{fig:c4p})
We construct the coupling $\Psi_{k}(x,y)$ of $D_{1}$ and $D_{2}$
using the following two step process. Let $\Psi_{j}$ be a coupling
of $D_{1}(j^{\prime})$ and $D_{2}(j^{\prime})$ which greedily maximises
the probability of assigning the same colour to site $j^{\prime}$
in each distribution. Then for each pair of colours $(c,c^{\prime})$
drawn from $\Psi_{j}$ we complete $\Psi_{k}(x,y)$ by letting it
be the coupling, greedily minimising Hamming distance, of the conditional
distributions $D_{1}\mid j^{\prime}=c$ and $D_{2}\mid j^{\prime}=c^{\prime}$
. 

\begin{figure}

\caption{Case \ref{c5}. \five.}

\label{fig:c5p} 

\begin{centering}\includegraphics{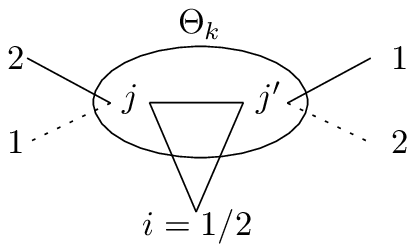} \par\end{centering}
\end{figure}

\begin{lem}
\label{lemma:tr2} Let $j$ and $j^{\prime}$ be the endpoints of
an edge $\Theta_{k}$and suppose that $\{ i,j\}\in E(G)$ and $\{ i,j^{\prime}\}\in E(G)$.
If $2$ is valid for $j$ in $D_{1}$, $1$ is valid for $j$ in $D_{2}$
and $1$ is not valid for $j^{\prime}$ in $D_{2}$ then \[
\rho_{i,j^{\prime}}^{k}\leq\frac{1}{q-\Delta+1}\textnormal{ and }\rho_{i,j}^{k}\leq\frac{1}{q-\Delta}.\]

\end{lem}
\begin{proof}
This is case \ref{c4} of the coupling. Note from Lemma \ref{lemma:d}
that $d_{1,j=2}=d_{2,j=1}$ and let $s=\sum_{c}d_{2,j=c}-d_{2,j=1}=\sum_{c\neq1}d_{2,j=c}$
denote the number of valid choices for $D_{2}$ other than choices
of the form $1c_{2}$. The number of valid choices for $D_{1}$ is
then $d_{1}=s+d_{1,j=2}+d_{1,j^{\prime}=2}$.

Since $1$ is not valid for $j^{\prime}$ in $D_{2}$ at least one
site other than $i$ on the boundary of $\Theta_{k}$ must be coloured
$1$ in $D_{1}$ (we say that some site $s$ on the boundary of $\Theta_{k}$
is coloured $c$ in $D_{1}$ if there exists a configuration with
positive measure in $D_{1}$ in which site $s$ is coloured $c$).
As there are no restrictions on the neighbourhood of $j$ each neighbour
of $j$ may be assigned a different colour in $D_{1}$. Hence we get
the following lower bounds on $d_{1,j=2}$ and $d_{1,j^{\prime}=2}$
\[
q-\Delta+1\leq d_{1,j=2}\textnormal{ and }q-\Delta\leq d_{1,j^{\prime}=2}.\]
 To lower bound $s$ observe that exactly $d_{1,j^{\prime}=2}$ colours
other than $1$ are valid for site $j$ in $D_{2}$ and let $J$ be
the set of colours, excluding $1$, that are valid for $j$ in $D_{2}$,
then \[
s=\sum_{c\in J}d_{2,j=c}\geq d_{1,j^{\prime}=c}\min_{j\in J}\{ d_{2,j=c}\}\geq d_{1,j^{\prime}=2}\left(d_{1,j=2}-1\right)\]
where we used the bound $d_{1,j=2}-1\leq d_{1,j=2}$ for $c\in J$
from Lemma \ref{lemma:triangle_d} $(i)$ since $1$ is not valid
for $j^{\prime}$ in $D_{2}$. 

We consider $\rho_{i,j^{\prime}}^{k}$ first. Suppose that a choice
of the form $c_{1}c_{2}$ is valid for $D_{2}$, in which case $c_{1}\neq2$
and $c_{2}\not\in\{1,2\}$ by the conditions of case \ref{c4} of
the coupling. Firstly if $c_{1}\neq1$ then $c_{1}c_{2}$ is also
valid for $D_{1}$ by Lemma \ref{lemma:d_same}. If $c_{1}=1$ then
the choice $2c_{2}$ is valid for $D_{1}$ by Lemma \ref{lemma:d}
and hence $d_{1}\geq d_{2}$. Note in particular that if a choice
$c_{1}c_{2}$ where $c_{2}\neq2$ is valid for $D_{1}$ then it is
also valid for $D_{2}$. Therefore, a different colour will only be
assigned to site $j^{\prime}$ in each distribution if $j^{\prime}$
is coloured 2 in $D_{1}$ during the first step of the coupling since
the Hamming distance at site $j^{\prime}$ is minimised greedily.
There are $d_{1,j^{\prime}=2}$ colourings assigning $2$ to $j^{\prime}$
in $D_{1}$ and hence \[
\rho_{i,j^{\prime}}^{k}\leq\frac{d_{1,j^{\prime}=2}}{d_{1,j=2}+d_{1,j^{\prime}=2}+s}\leq\frac{d_{1,j^{\prime}=2}}{d_{1,j=2}+\left(1+d_{1,j^{\prime}=2}\right)}\leq\frac{1}{d_{1,j=2}}\leq\frac{1}{q-\Delta+1}\]
 using the lower bounds on $s$ and $d_{1,j=2}$.

Now consider $\rho_{i,j}^{k}$. Suppose that $_{1}^{\prime},c_{2}^{\prime})$
is the pair of colours drawn for site $j^{\prime}$ in the first step
of the coupling. The second step of $\Psi_{k}(x,y)$ then couples
the conditional distributions $D_{1}\mid j^{\prime}=c_{1}^{\prime}$
and $D_{2}\mid j^{\prime}=c_{2}^{\prime}$ greedily to minimise Hamming
distance. First suppose that $c_{1}^{\prime}\neq c_{2}^{\prime}$.
It was pointed out in the analysis above that if $c_{1}^{\prime}\neq c_{2}^{\prime}$
then $c_{1}^{\prime}=2$ and the resulting configuration is shown
in Figure \ref{fig:intermediate-step}. %
\begin{figure}

\caption{The pair of configurations after the colour of site $j^{\prime}$
has been assigned during the first step of the coupling.}

\begin{centering}\label{fig:intermediate-step}\includegraphics{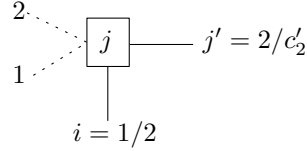}\par\end{centering}
\end{figure}
 We make the following observations about the resulting conditional
distributions $D_{1}\mid j^{\prime}=2$ and $D_{2}\mid j^{\prime}=c_{2}^{\prime}$.
\begin{itemize}
\item The colour $2$ is \emph{not} valid for $j$ in either distribution
$D_{1}\mid j^{\prime}=2$ or $D_{2}\mid j^{\prime}=c_{2}^{\prime}$.
\item The colour $1$ is \emph{not} valid for $j$ in $D_{1}\mid j^{\prime}=2$
but \emph{could} be valid for $j$ in $D_{2}\mid j^{\prime}=c_{2}^{\prime}$.
\item The colour $c_{2}^{\prime}$ \emph{could} be valid for $j$ in $D_{1}\mid j^{\prime}=2$
but is \emph{not} valid for $j$ in $D_{2}\mid j^{\prime}=c_{2}^{\prime}$.
\item For each $c\in C\setminus\{1,2,c_{2}^{\prime}\}$ the colour $c$
is valid for $j$ in $D_{1}\mid j^{\prime}=2$ \emph{if and only if}
$c$ is valid for $j$ in $D_{2}\mid j^{\prime}=c_{2}^{\prime}$.
\end{itemize}
These observations show that this case is a single-site disagreement
sub problem and that there must be at least $(q-3)-(\Delta-2)=q-\Delta-1$
colours that are valid for $j$ in both conditional distributions
since $j$ has at most $\Delta-2$ neighbours other than $i$ and
$j^{\prime}$. Also, there is at most one colour which is valid for
$j$ in one distribution but not in the other and since the coupling
greedily maximises Hamming distance this implies \[
\Prob_{(x^{\prime},y^{\prime})\in\Psi_{k}(x,y)}(x_{j}^{\prime}\neq y_{j}^{\prime}\mid x_{j^{\prime}}^{\prime}\neq y_{j^{\prime}}^{\prime})\leq\frac{1}{q-\Delta}.\]

Now suppose that the same colour $c$, say, is drawn for site $j^{\prime}$
in both distributions during the first step of the coupling. Then
the only site adjacent to $i$ that is coloured differently in the
conditional distributions $D_{1}\mid j^{\prime}=c$ and $D_{2}\mid j^{\prime}=c$
is site $i$, so using a similar reasoning to above we find \[
\Prob_{(x^{\prime},y^{\prime})\in\Psi_{k}(x,y)}(x_{j}^{\prime}\neq y_{j}^{\prime}\mid x_{j^{\prime}}^{\prime}=y_{j^{\prime}}^{\prime})\leq\frac{1}{q-\Delta}\]
and thus \begin{align*}
\rho_{i,j}^{k} & =\max_{(x,y)\in S_{i}}\left\{ \Prob_{(x^{\prime},y^{\prime})\in\Psi_{k}(x,y)}(x_{j}^{\prime}\neq y_{j}^{\prime})\right\} \\
 & =\max_{(x,y)\in S_{i}}\{\Prob_{(x^{\prime},y^{\prime})\in\Psi_{k}(x,y)}(x_{j}^{\prime}\neq y_{j}^{\prime}\mid x_{j^{\prime}}^{\prime}\neq y_{j^{\prime}}^{\prime})\Prob_{(x^{\prime},y^{\prime})\in\Psi_{k}(x,y)}(x_{j^{\prime}}^{\prime}\neq y_{j^{\prime}}^{\prime})\\
 & \quad+\Prob_{(x^{\prime},y^{\prime})\in\Psi_{k}(x,y)}(x_{j}^{\prime}\neq y_{j}^{\prime}\mid x_{j^{\prime}}^{\prime}=y_{j^{\prime}}^{\prime})\Prob_{(x^{\prime},y^{\prime})\in\Psi_{k}(x,y)}(x_{j^{\prime}}^{\prime}=y_{j^{\prime}}^{\prime})\}\\
 & \leq\max_{(x,y)\in S_{i}}\left\{ \frac{1}{q-\Delta}\Prob_{(x^{\prime},y^{\prime})\in\Psi_{k}(x,y)}(x_{j^{\prime}}^{\prime}\neq y_{j^{\prime}}^{\prime})+\frac{1}{q-\Delta}\Prob_{(x^{\prime},y^{\prime})\in\Psi_{k}(x,y)}(x_{j^{\prime}}^{\prime}=y_{j^{\prime}}^{\prime})\right\} \\
 & =\frac{1}{q-\Delta}\max_{(x,y)\in S_{i}}\left\{ \Prob_{(x^{\prime},y^{\prime})\in\Psi_{k}(x,y)}(x_{j^{\prime}}^{\prime}\neq y_{j^{\prime}}^{\prime})+\Prob_{(x^{\prime},y^{\prime})\in\Psi_{k}(x,y)}(x_{j^{\prime}}^{\prime}=y_{j^{\prime}}^{\prime})\right\} =\frac{1}{q-\Delta}\end{align*}
 which completes the proof.
\end{proof}
\textbf{Details of case \ref{c5}.} (Repeated in Figure \ref{fig:c5p})
 We assume without loss of generality that $d_{1}\geq d_{2}$ and
construct the coupling $\Psi_{k}(x,y)$ of $D_{1}$ and $D_{2}$ as
follows. For each valid choice of the form $c_{1}c_{2}$ for $D_{1}$
where $c_{1}\neq1$ and $c_{2}\neq2$ Lemma \ref{lemma:d_same} guarantees
that $c_{1}c_{2}$ is also valid for $D_{2}$ so we construct $\Psi_{k}(x,y)$
such that \begin{align*}
\Prob_{(x^{\prime},y^{\prime})\in\Psi_{k}(x,y)}(x^{\prime}=y^{\prime}=c_{1}c_{2})= & \frac{1}{d_{1}}.\end{align*}
This leaves the set $Z_{1}=\{ c_{1}2\mid c_{1}2\textnormal{ valid in }D_{1}\}$
of valid choices for $D_{1}$ and $Z_{2}=\{1c_{2}\mid1c_{2}\textnormal{ valid in }D_{2}\}\subseteq D_{2}$
for $D_{2}$. Observe that $z_{1}\geq z_{2}$ where $z_{1}$ and $z_{2}$
denote the size of $Z_{1}$ and $Z_{2}$ respectively. Let $Z_{1}(t)$
denote the $t$-th element of $Z_{1}$ and similarly for $Z_{2}$.
Then for $1\leq t\leq z_{2}$ let\[
\Prob_{(x^{\prime},y^{\prime})\in\Psi_{k}(x,y)}(x^{\prime}=Z_{1}(t),y^{\prime}=Z_{2}(t))=\frac{1}{d_{1}}\]
and for each pair $z_{2}+1\leq t\leq z_{2}$ and $h\in D_{2}$ let\[
\Prob_{(x^{\prime},y^{\prime})\in\Psi_{k}(x,y)}(x^{\prime}=Z_{1}(t),y^{\prime}=h)=\frac{1}{d_{1}d_{2}}.\]
It is easy to verify that each valid colouring has the correct weight
in $\Psi_{k}(x,y)$ so this completes the coupling. 

\begin{figure}

\caption{Case \ref{c6}. \six.}

\label{fig:c6p} 

\begin{centering}\includegraphics{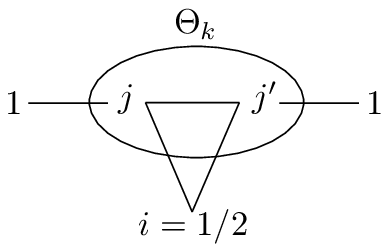} \par\end{centering}
\end{figure}

\begin{lem}
\label{lemma:tr3} Let $j$ and $j^{\prime}$ be the endpoints of
an edge $\Theta_{k}$ and suppose that $\{ i,j\}\in E(G)$ and $\{ i,j^{\prime}\}\in E(G)$.
If $1$ is valid for $j$ in $D_{2}$, $1$ is not valid for $j^{\prime}$
in $D_{2}$ and $2$ is not valid for $j$ in $D_{1}$ then \[
\rho_{i,j}^{k}\leq\rho_{i,j^{\prime}}^{k}\leq\frac{1}{q-\Delta}.\]

\end{lem}
\begin{proof}
This is case \ref{c5} of the coupling. Let $s=\sum_{c}d_{2,j=c}-d_{2,j=1}$
be the number of valid choices for $D_{2}$ other than choices of
the form $1c_{2}$. Observe that $d_{2}=s+d_{1,j^{\prime}=2}$ and
note that $d_{1,j^{\prime}=2}\geq d_{2,j=1}$ since we have assumed
$d_{1}\geq d_{2}$ in the construction of the coupling. At least one
neighbour, other than $i$, of $j^{\prime}$ on the boundary of $\Theta_{k}$
is coloured 1 in $D_{1}$ and we get the following lower-bound on
$d_{2,j=1}$ since all other neighbours of $j^{\prime}$ may be assigned
a different colour \[
q-\Delta+1\leq d_{2,j=1}.\]
 We bound $s$ using the same argument as in the proof of Lemma \ref{lemma:tr2}
and get \[
d_{1,j^{\prime}=2}(d_{2,j=1}-1)\leq s.\]

Since $2$ is not valid for $j$ in $D_{1}$ the first $d_{2,j=1}$
choices of the form $c_{1}2$ for $D_{1}$ are matched with some choice
of the form $1c_{1}$ for $D_{2}$ with probability $1/d_{1}$ resulting
in a different colour being assigned to both sites $j$ and $j^{\prime}$
in each distribution. Each of the $d_{1,j^{\prime}=2}-d_{2,j=1}$
remaining valid choices for $D_{1}$ is matched with each valid choice
for $D_{2}$ with probability $\frac{1}{d_{1}d_{2}}$ resulting in
a disagreement at $j^{\prime}$ (since $2$ is not valid for $j^{\prime}$
in $D_{2}$) and potentially also at at $j$ so $\rho_{i,j}^{k}\leq\rho_{i,j^{\prime}}^{k}$.
Hence the probability of making a choice of the form $c_{1}2$ for
$D_{1}$\[
\Prob_{(x^{\prime},y^{\prime})\in\Psi_{k}(x,y)}(x_{j^{\prime}}^{\prime}=2)=\frac{d_{2,j=1}}{d_{1}}+\sum_{h\in D_{2}}\frac{d_{1,j^{\prime}=2}-d_{2,j=1}}{d_{1}d_{2}}=\frac{d_{1,j^{\prime}=2}}{d_{1}}\]
is an upper bound on the disagreement probabilities at both sites
$j$ and $j^{\prime}$. Using the lower bounds on $s$ and $d_{2,j=1}$
we have \[
\rho_{i,j}^{k}\leq\rho_{i,j^{\prime}}^{k}\leq\frac{d_{1,j^{\prime}=2}}{d_{1}}=\frac{d_{1,j^{\prime}=2}}{d_{1,j^{\prime}=2}+s}\leq\frac{d_{1,j^{\prime}=2}}{d_{1,j^{\prime}=2}+(d_{2,j=1}-1)d_{1,j^{\prime}=2}}\leq\frac{1}{q-\Delta}\]
which completes the proof.
\end{proof}
\textbf{Details of case \ref{c6}.} (Repeated in Figure \ref{fig:c6p})
 First observe that $1$ is not valid for both $j$ and $j^{\prime}$
in either distribution $D_{1}$ or $D_{2}$ so $d_{1}=d_{2}+d_{1,j=2}+d_{1,j^{\prime}=2}\geq d_{2}$
by Lemma \ref{lemma:d_same}, since any choice valid for $D_{2}$
does not assign colour $2$ to any site in $\Theta_{k}$. Let $Z_{1}$
and $Z_{2}$ be the sets of colourings valid for $D_{1}$ and $D_{2}$
respectively. We define the following mutually exclusive subsets of
$Z_{1}$. $Z_{j}=\{2c_{2}\mid2c_{2}\in Z_{1}\}$, $Z_{j^{\prime}}=\{ c_{1}2\mid c_{1}2\in Z_{1}\}$
and $Z=Z_{1}\setminus(Z_{j}\cup Z_{j^{\prime}})=Z_{2}$. By construction,
the union of these three subsets is $Z_{1}$ and note that the size
of $Z_{j}$ is $d_{1,j=2}$, the size of $Z_{j^{\prime}}$ is $d_{1,j^{\prime}=2}$
and the size of $Z$ is $d_{2}$. 

First we consider choices from $Z$ for $D_{1}$. For each choice
$h\in Z$ we have $h\in Z_{2}$ by construction of $Z$ and so we
use the identity coupling and let

\[
\Prob_{(x^{\prime},y^{\prime})\in\Psi_{k}(x,y)}(x^{\prime}=y^{\prime}=h)=\frac{1}{d_{1}}.\]
We let the remainder of the coupling minimise Hamming distance. First
consider the choices for $D_{1}$ in $Z_{j}$. We construct $\Psi_{k}(x,y)$
such that it minimises Hamming distance and assigns probability $1/d_{1}$
to each choice for $D_{1}$ in $Z_{j}$ whilst ensuring that for each
choice $g\in Z_{2}$ for $D_{2}$ \[
\sum_{h\in Z_{j}}\Prob_{(x^{\prime},y^{\prime})\in\Psi_{k}(x,y)}(x^{\prime}=h,y^{\prime}=g)=\frac{d_{1,j=2}}{d_{1}d_{2}}.\]
Similarly we assign probability $1/d_{1}$ to each choice for $D_{1}$
in $Z_{j^{\prime}}$ whilst also requiring that for each choice $g\in Z_{2}$
for $D_{2}$\[
\sum_{h\in Z_{j^{\prime}}}\Prob_{(x^{\prime},y^{\prime})\in\Psi_{k}(x,y)}(x^{\prime}=h,y^{\prime}=g)=\frac{d_{1,j^{\prime}=2}}{d_{1}d_{2}}.\]

To see that this ensures that the coupling is fair observe that each
choice $h\in Z_{1}$ receives weight $1/d_{1}$ and each choice $g\in Z_{2}$
weight \[
\frac{1}{d_{1}}+\frac{d_{1,j=2}}{d_{1}d_{2}}+\frac{d_{1,j^{\prime}=2}}{d_{1}d_{2}}=\frac{d_{2}+d_{1,j=2}+d_{1,j^{\prime}=2}}{d_{1}d_{2}}=\frac{1}{d_{2}}\]
since $d_{2}+d_{1,j=2}+d_{1,j^{\prime}=2}=d_{1}$.

\begin{rem*}
Note that a coupling satisfying these requirements always exists.
We will not give the detailed construction of $\Psi_{k}(x,y)$ here,
but in the subsequent proof we will consider three cases. In the first
two cases \emph{any} coupling minimising Hamming distance will be
sufficient to establish the required bounds on the influence of $i$
on $j$. In the final case we will need a detailed construction of
the coupling and so will provide it together with the proof for ease
of reference.
\end{rem*}
\begin{lem}
\label{lemma:tr4} Let $j$ and $j^{\prime}$ be the endpoints of
an edge $\Theta_{k}$ and suppose that $\{ i,j\}\in E(G)$ and $\{ i,j^{\prime}\}\in E(G)$.
If $1$ is not valid for $j$ in $D_{2}$ and $1$ is not valid for
$j^{\prime}$ in $D_{2}$ then \[
\rho_{i,j}^{k}\leq\frac{1}{q-\Delta+1}+\frac{1}{(q-\Delta+1)^{2}}\textnormal{ and }\rho_{i,j^{\prime}}^{k}\leq\frac{1}{q-\Delta+1}+\frac{1}{(q-\Delta+1)^{2}}.\]

\end{lem}
\begin{proof}
This is case \ref{c6} of the coupling. We consider three separate
cases. Firstly suppose that $2$ is not valid for both $j$ and $j^{\prime}$
in $D_{1}$. Then the only valid choices for $D_{1}$ are of the form
$c_{1}c_{2}$ where $c_{1},c_{2}\in C\setminus\{1,2\}$ and each such
choice is also valid in $D_{2}$ as observed in the construction of
the coupling. The same colouring is selected for each distribution
and hence \[
\rho_{i,j}^{k}=0\textnormal{ and }\rho_{i,j^{\prime}}^{k}=0.\]

Next suppose that exactly one site in $\Theta_{k}$, $j^{\prime}$
say, is adjacent to some site coloured $2$ in $D_{1}$. As in the
previous case, each choice that is valid in both $D_{1}$ and $D_{2}$
is matched using the identity matching and does not cause a discrepancy
at any site. However if a choice of the form $2c$ is made for $D_{1}$
then site $j$ will be coloured differently in each colouring drawn
from $\Psi_{k}(x,y)$ and the colour at site $j^{\prime}$ may also
be different so $\rho_{i,j^{\prime}}^{k}\leq\rho_{i,j}^{k}$. Since
all choices of the form $c2$ are not valid for $D_{1}$, making a
choice of the form $2c$ for $D_{1}$ is the only way to create a
disagreement at any site in the coupling and so \[
\rho_{i,j^{\prime}}^{k}\leq\rho_{i,j}^{k}\leq\frac{d_{1,j=2}}{d_{1}}\]
since $d_{1,j=2}$ is the number of valid choices for $D_{1}$ of
the form $2c$. We need to establish a lower bound of $d_{1}$ and
observe that, for $c$ valid for $j$ in $D_{1}$, $d_{1,j=2}-1\leq d_{1,j=c}$
by Lemma \ref{lemma:triangle_d} (ii) since $2$ is not valid for
$j^{\prime}$ in $D_{1}$. Let $v$ be the number of colours that
are valid for site $j$ in $D_{1}$. Then $v$ is lower bounded by
$q-\Delta+2\leq v$ since at least two of the sites (including $i$)
adjacent to $j$ on the boundary of $\Theta_{k}$ are coloured $1$
in $D_{1}$. Also, since at least one site (other than $j$ and $i$)
adjacent to $j^{\prime}$ is coloured $1$ and another is coloured
$2$ in $D_{1}$, we have $q-\Delta+2\leq d_{1,j=2}$. Using the lower
bounds on $v$ and $d_{1,j=c}$ we have, letting $J$ denote the set
of colours other than $2$ that are valid for $j$ in $D_{1}$,\begin{align*}
d_{1}=\sum_{c}d_{1,j=c} & =d_{1,j=2}+\sum_{c\in J}d_{1,j=c}\geq d_{1,j=2}+\sum_{c\in J}(d_{1,j=2}-1)\\
 & \geq(v-1)(d_{1,j=2}-1)+d_{1,j=2}\geq(q-\Delta+2)d_{1,j=2}-(q-\Delta+1)\end{align*}
and hence using the lower bound on $d_{1,j=2}$ \begin{align*}
\frac{1}{\rho_{i,j}^{k}} & \geq\frac{(q-\Delta+2)d_{1,j=2}-(q-\Delta+1)}{d_{1,j=2}}\geq q-\Delta+2-\frac{q-\Delta+1}{q-\Delta+2}>q-\Delta+1\end{align*}
which gives the bounds required by the statement of the lemma.

Finally consider the case when the colour $2$ is valid for both $j$
and $j^{\prime}$ in $D_{1}$. In this case we will provide details
of the construction of $\Psi_{k}(x,y)$ when required. We begin by
establishing some required bounds. Since $1$ is not valid for $j^{\prime}$
in $D_{2}$ at least two neighbours of $j^{\prime}$ (including $i$)
must be coloured $1$ in $D_{1}$ and the same applies to the neighbourhood
of $j$, so we get the following lower bounds on $d_{1,j=2}$ and
$d_{1,j^{\prime}=2}$ \[
q-\Delta+1\leq d_{1,j=2}\textnormal{ and }q-\Delta+1\leq d_{1,j^{\prime}=2}.\]
 We also require bounds on $d_{2,j=c}$ and $d_{2,j^{\prime}=c}$
for other colours $c$. Suppose that the choice $cc^{\prime}$ is
valid in $D_{2}$ then, since $c,c^{\prime}\in C\setminus\{1,2\}$,
$cc^{\prime}$ is also valid for $D_{1}$ by Lemma \ref{lemma:d_same}.
Furthermore, the choice $c2$ is valid in $D_{1}$ (but not $D_{2}$)
so $d_{1,j=c}-1=d_{2,j=c}$. Lemma \ref{lemma:triangle_d} (ii) guarantees
that $d_{1,j=2}\leq d_{1,j=c}\leq d_{1,j=2}+1$ so \[
d_{1,j=2}-1\leq d_{2,j=c}\leq d_{1,j=2}\]
for any $c$ valid for $j$ in $D_{1}$ and a similar argument gives
the bound \[
d_{1,j^{\prime}=2}-1\leq d_{2,j^{\prime}=c}\leq d_{1,j^{\prime}=c}\]
for any colour $c$ valid for $j^{\prime}$ in $D_{2}$. Observe that
exactly $d_{1,j^{\prime}=2}$ colours must be valid for site $j$
in $D_{2}$ so using the stated bounds on $d_{2,j=c}$ we have the
following bounds on $d_{2}$ \[
d_{1,j^{\prime}=2}(d_{1,j=2}-1)\leq d_{2}\leq d_{1,j^{\prime}=2}d_{1,j=2}.\]

We bound the probability of disagreements at sites $j$ and $j^{\prime}$
from choices made for $D_{1}$. From the coupling we again note that
if a choice $c_{1}c_{2}$ where $c_{1}\neq2$ and $c_{2}\neq2$ is
made for $D_{1}$ then there will be no disagreements at any site
in $\Theta_{k}$.

Consider making a valid choice of the form $2c$ for $D_{1}$. Firstly,
such a choice for $D_{1}$ will cause site $j$ to be coloured differently
in any pair of colourings drawn from the coupling since $2$ is not
valid for $j$ in $D_{2}$. We construct $\Psi_{k}(x,y)$ such that
the choice $2c$ for $D_{1}$ is matched with a choice of the form
$c^{\prime}c$ for $D_{2}$ as long as such a choice that has not
exceeded it aggregated probability exists. Let $J$ denote the set
of choices of the form $c^{\prime}c$ that are valid for $D_{2}$
and note that the size of $J$ is $d_{2,j^{\prime}=c}$. The total
aggregated weight of all choices of the form $c^{\prime}c$ for $D_{2}$
is\[
\sum_{g\in J}\sum_{h\in Z_{1}}\Prob_{(x^{\prime},y^{\prime})\in\Psi_{k}(x,y)}(x^{\prime}=h,y^{\prime}=g)=\sum_{g\in J}\frac{d_{1,j=2}}{d_{1}d_{2}}=\frac{d_{2,j^{\prime}=c}d_{1,j=2}}{d_{1}d_{2}}\]
so as long as \[
\frac{1}{d_{1}}\leq\frac{d_{2,j^{\prime}=c}d_{1,j=2}}{d_{1}d_{2}}\]
there is enough probability available in $Z_{2}$ to match all the
weight of the choice $2c$ for $D_{1}$ with a choice of the form
$c^{\prime}c$ for $D_{2}$ and hence assigning the same colour, $c$,
to site $j^{\prime}$ in any pair of colourings drawn from the coupling.
If there is not enough unassigned weight available in $Z_{2}$ then
the coupling will match at much probability as possible, $\frac{d_{2,j^{\prime}=c}d_{1,j=2}}{d_{1}d_{2}}$,
with choices of the form $c^{\prime}c$ for $Z_{2}$ but the remaining
probability will be matched with choices not assigning colour $c$
to site $j^{\prime}$ in $Z_{2}$. Hence we obtain the following probabilities
conditioned on making a choice of the form $2c$ for $D_{1}$. \[
\Prob_{(x^{\prime},y^{\prime})\in\Psi_{k}(x,y)}(x_{j}^{\prime}\neq y_{j}^{\prime}\mid x^{\prime}=2c)=1\]
 and \begin{align*}
\Prob_{(x^{\prime},y^{\prime})\in\Psi_{k}(x,y)}(x_{j^{\prime}}^{\prime}\neq y_{j^{\prime}}^{\prime}\mid x^{\prime}=2c) & \leq\max\left(0,1-\frac{d_{2,j^{\prime}=c}d_{1,j=2}}{d_{2}}\right)\\
 & \leq\max\left(0,1-\frac{(d_{1,j^{\prime}=2}-1)d_{1,j=2}}{d_{1,j=2}d_{1,j^{\prime}=2}}\right)\\
 & \leq\frac{1}{d_{1,j^{\prime}=2}}\end{align*}
 using the bounds on $d_{2}$ and $d_{1,j^{\prime}=c}$. Lastly observe
that there are $d_{1,j=2}$ valid choices for $D_{1}$ of the form
$2c$ so \[
\sum_{c}\Prob_{(x^{\prime},y^{\prime})\in\Psi_{k}(x,y)}(x^{\prime}=2c)=\frac{d_{1,j=2}}{d_{1}}=\frac{d_{1,j=2}}{d_{1,j=2}+d_{1,j^{\prime}=2}+d_{2}}.\]

Also consider making a valid choice of the form $c2$ for $D_{1}$.
This case is symmetric to the construction above, but we include it
for completeness. A choice of the form $c2$ for $D_{1}$ will cause
site $j^{\prime}$ to be coloured differently in any pair of colourings
drawn from the coupling since $2$ is not valid for $j^{\prime}$
in $D_{2}$. We hence construct $\Psi_{k}(x,y)$ such that it matches
the choice $c2$ for $D_{1}$ with a choice of the form $cc^{\prime}$
for $D_{2}$ as long as such a choice that has not exceeded it aggregated
probability exists. Let $J^{\prime}$ denote the set of choices of
the form $cc^{\prime}$ that are valid for $D_{2}$ and note that
the size of $J^{\prime}$ is $d_{2,j=c}$. The total aggregated weight
of all choices of the form $cc^{\prime}$ for $D_{2}$ is\[
\sum_{g\in J^{\prime}}\sum_{h\in Z_{1}}\Prob_{(x^{\prime},y^{\prime})\in\Psi_{k}(x,y)}(x^{\prime}=h,y^{\prime}=g)=\sum_{g\in J}\frac{d_{1,j^{\prime}=2}}{d_{1}d_{2}}=\frac{d_{2,j=c}d_{1,j^{\prime}=2}}{d_{1}d_{2}}\]
so as long as \[
\frac{1}{d_{1}}\leq\frac{d_{2,j=c}d_{1,j^{\prime}=2}}{d_{1}d_{2}}\]
there is enough weight available in $Z_{2}$ to match all the weight
of the choice $c2$ for $D_{1}$ with a choice of the form $cc^{\prime}$
for $D_{2}$ and hence assigning the same colour, $c$, to site $j$
in any pair of colourings drawn from the coupling. If there is not
enough unassigned weight available in $Z_{2}$ then the coupling will
match at much weight as possible, $\frac{d_{2,j=c}d_{1,j^{\prime}=2}}{d_{1}d_{2}}$,
with choices of the form $cc^{\prime}$ for $Z_{2}$ but the remaining
weight will be matched with choices not assigning colour $c$ to site
$j$ in $Z_{2}$. Hence we obtain the following probabilities conditioned
on making a choice of the form $c2$ for $D_{1}$ 

\begin{align*}
\Prob_{(x^{\prime},y^{\prime})\in\Psi_{k}(x,y)}(x_{j}^{\prime}\neq y_{j}^{\prime}\mid x^{\prime}=c2) & \leq\max\left(0,1-\frac{d_{2,j=c}d_{1,j^{\prime}=2}}{d_{2}}\right)\\
 & \leq\max\left(0,1-\frac{d_{1,j^{\prime}=2}(d_{1,j=2}-1)}{d_{1,j=2}d_{1,j^{\prime}=2}}\right)\\
 & \leq\frac{1}{d_{1,j=2}}\end{align*}
 using the bounds on $d_{2}$ and $d_{1,j=c}$, and as before we also
have \[
\Prob_{(x^{\prime},y^{\prime})\in\Psi_{k}(x,y)}(x_{j^{\prime}}^{\prime}\neq y_{j^{\prime}}^{\prime}\mid x^{\prime}=c2)=1\mbox{ and }\sum_{c}\Prob_{(x^{\prime},y^{\prime})\in\Psi_{k}(x,y)}(x^{\prime}=c2)=\frac{d_{1,j^{\prime}=2}}{d_{1,j=2}+d_{1,j^{\prime}=2}+d_{2}}.\]
Using the conditional probabilities and the bounds on $d_{2}$, $d_{1,j=2}$
and $d_{1,j^{\prime}=2}$ we find \begin{align*}
\rho_{i,j}^{k} & =\max_{(x,y)\in S_{i}}\left\{ \Prob_{(x^{\prime},y^{\prime})\in\Psi_{k}(x,y)}(x_{j}^{\prime}\neq y_{j}^{\prime})\right\} \\
 & =\max_{(x,y)\in S_{i}}\Big\{\sum_{c}\Prob_{(x^{\prime},y^{\prime})\in\Psi_{k}(x,y)}(x_{j}^{\prime}\neq y_{j}^{\prime}\mid x^{\prime}=2c)\Prob_{(x^{\prime},y^{\prime})\in\Psi_{k}(x,y)}(x^{\prime}=2c)\\
 & \quad+\sum_{c}\Prob_{(x^{\prime},y^{\prime})\in\Psi_{k}(x,y)}(x_{j}^{\prime}\neq y_{j}^{\prime}\mid x^{\prime}=c2)\Prob_{(x^{\prime},y^{\prime})\in\Psi_{k}(x,y)}(x^{\prime}=c2)\Big\}\\
 & \leq\max_{(x,y)\in S_{i}}\left\{ \sum_{c}\left[\Prob_{(x^{\prime},y^{\prime})\in\Psi_{k}(x,y)}(x^{\prime}=2c)+\Prob_{(x^{\prime},y^{\prime})\in\Psi_{k}(x,y)}(x^{\prime}=c2)\frac{1}{d_{1,j=2}}\right]\right\} \\
 & \leq\max_{(x,y)\in S_{i}}\left\{ \frac{d_{1,j=2}}{d_{1,j=2}+d_{1,j^{\prime}=2}+d_{2}}+\sum_{c}\Prob_{(x^{\prime},y^{\prime})\in\Psi_{k}(x,y)}(x^{\prime}=c2)\frac{1}{d_{1,j=2}}\right\} \\
 & \leq\max_{(x,y)\in S_{i}}\left\{ \frac{d_{1,j=2}}{d_{1,j=2}(1+d_{1,j^{\prime}=2})}+\frac{d_{1,j^{\prime}=2}}{d_{1,j=2}+d_{1,j^{\prime}=2}+d_{2}}\frac{1}{d_{1,j=2}}\right\} \\
 & \leq\max_{(x,y)\in S_{i}}\left\{ \frac{1}{1+d_{1,j^{\prime}=2}}+\frac{d_{1,j^{\prime}=2}}{d_{1,j=2}(1+d_{1,j^{\prime}=2})}\frac{1}{d_{1,j=2}}\right\} \\
 & \leq\frac{1}{q-\Delta+2}+\frac{1}{(q-\Delta+1)^{2}}\end{align*}
 and again by symmetry \begin{align*}
\rho_{i,j^{\prime}}^{k} & =\max_{(x,y)\in S_{i}}\left\{ \Prob_{(x^{\prime},y^{\prime})\in\Psi_{k}(x,y)}(x_{j^{\prime}}^{\prime}\neq y_{j^{\prime}}^{\prime}\right\} \\
 & =\max_{(x,y)\in S_{i}}\Big\{\sum_{c}\Prob_{(x^{\prime},y^{\prime})\in\Psi_{k}(x,y)}(x_{j^{\prime}}^{\prime}\neq y_{j^{\prime}}^{\prime}\mid x^{\prime}=2c)\Prob_{(x^{\prime},y^{\prime})\in\Psi_{k}(x,y)}(x^{\prime}=2c)\\
 & \quad+\sum_{c}\Prob_{(x^{\prime},y^{\prime})\in\Psi_{k}(x,y)}(x_{j^{\prime}}^{\prime}\neq y_{j^{\prime}}^{\prime}\mid x^{\prime}=c2)\Prob_{(x^{\prime},y^{\prime})\in\Psi_{k}(x,y)}(x^{\prime}=c2)\Big\}\\
 & \leq\max_{(x,y)\in S_{i}}\left\{ \sum_{c}\left[\Prob_{(x^{\prime},y^{\prime})\in\Psi_{k}(x,y)}(x^{\prime}=2c)+\Prob_{(x^{\prime},y^{\prime})\in\Psi_{k}(x,y)}(x^{\prime}=c2)\frac{1}{d_{1,j^{\prime}=2}}\right]\right\} \\
 & \leq\max_{(x,y)\in S_{i}}\left\{ \frac{d_{1,j^{\prime}=2}}{d_{1,j=2}+d_{1,j^{\prime}=2}+d_{2}}+\sum_{c}\Prob_{(x^{\prime},y^{\prime})\in\Psi_{k}(x,y)}(x^{\prime}=c2)\frac{1}{d_{1,j^{\prime}=2}}\right\} \\
 & \leq\max_{(x,y)\in S_{i}}\left\{ \frac{d_{1,j^{\prime}=2}}{d_{1,j=2}(1+d_{1,j^{\prime}=2})}+\frac{d_{1,j=2}}{d_{1,j=2}+d_{1,j^{\prime}=2}+d_{2}}\frac{1}{d_{1,j^{\prime}=2}}\right\} \\
 & \leq\max_{(x,y)\in S_{i}}\left\{ \frac{1}{1+d_{1,j^{=}2}}+\frac{d_{1,j=2}}{d_{1,j=2}(1+d_{1,j^{\prime}=2})}\frac{1}{d_{1,j^{\prime}=2}}\right\} \\
 & \leq\frac{1}{q-\Delta+1}+\frac{1}{(q-\Delta+1)^{2}}\end{align*}
 which implies the statement of the lemma.
\end{proof}
This completes the cases of the coupling and we combine the obtained
bounds on $\rho_{i,j}^{k}$ and $\rho_{i,j^{\prime}}^{k}$ in the
following corollary of Lemmas \ref{lemma:tr1}, \ref{lemma:tr2},
\ref{lemma:tr3} and \ref{lemma:tr4} which we use in establishing
the mixing time of $\escan$.

\begin{cor}
\label{cor:triangle} Let $j$ and $j^{\prime}$ be the endpoints
of an edge $\Theta_{k}$. If $\{ i,j\}\in E(G)$ and $\{ i,j^{\prime}\}\in E(G)$
then \[
\rho_{i,j}^{k}\leq\frac{1}{q-\Delta}+\frac{1}{(q-\Delta)^{2}}\textnormal{ and }\rho_{i,j^{\prime}}^{k}\leq\frac{1}{q-\Delta}+\frac{1}{(q-\Delta)^{2}}.\]

\end{cor}
\begin{rem*}
Note that the bound in Corollary \ref{cor:triangle} is never tight.
This bound could be improved, however this would only allow us to
beat the $2\Delta$ bound for special graphs since the bounds in Lemma
\ref{lemma:gen2delta_notriag} are tight.
\end{rem*}
We are now ready to present a proof of Theorem \ref{thm:gen2delta}.
\emph{}\textbf{\\[1 em]Theorem \ref{thm:gen2delta}.} \emph{Let $G$
be a graph with maximum vertex-degree $\Delta$. If $q\geq2\Delta$
then} \[
\Mix(\escan,\epsilon)\leq\Delta^{2}\log(n\epsilon^{-1}).\]

\begin{proof}
Let $j$ and $j^{\prime}$ be the endpoints of an edge represented
by a block $\Theta_{k}$. Let $\alpha_{j}=\sum_{i}\rho_{i,j}^{k}$
be the influence on site $j$ and $\alpha_{j^{\prime}}=\sum_{i}\rho_{i,j^{\prime}}^{k}$
then influence on $j^{\prime}$. Then $\alpha=\max(\alpha_{j},\alpha_{j^{\prime}})$.
Suppose that $\Theta_{k}$ is adjacent to $t$ triangles, that is
there are $t$ sites $i_{1},\dots,i_{t}$ such that $\{ i,j\}\in E(G)$
and $\{ i,j^{\prime}\}\in E(G)$ for each $i\in\{ i_{1},\dots,i_{t}\}$.
Note that $0\leq t\leq\Delta-1$. There are at most $\Delta-1-t$
sites adjacent to $j$ that are not adjacent to $j^{\prime}$ and
at most $\Delta-1-t$ sites adjacent to $j^{\prime}$ that are not
adjacent to $j$. From Lemma \ref{lemma:gen2delta_notriag} a site
adjacent only to $j$ will emit an influence of at most $\frac{1}{q-\Delta}$
on site $j$ and Lemma \ref{lemma:gen2delta_notriag} also guarantees
that a site only adjacent to $j^{\prime}$ can emit an influence at
most $\frac{1}{(q-\Delta)^{2}}$ on site $j$. Corollary \ref{cor:triangle}
says that a site adjacent to both $j$ and $j^{\prime}$ can emit
an influence of at most $\frac{1}{q-\Delta}+\frac{1}{(q-\Delta)^{2}}$
on site $j$ and hence \begin{align*}
\alpha_{j} & \leq t\left(\frac{1}{q-\Delta}+\frac{1}{(q-\Delta)^{2}}\right)+(\Delta-1-t)\left(\frac{1}{q-\Delta}\right)+(\Delta-1-t)\left(\frac{1}{(q-\Delta)^{2}}\right)\\
 & =\frac{\Delta-1}{q-\Delta}+\frac{\Delta-1}{(q-\Delta)^{2}}\end{align*}
and similarly by considering the influence on site $j^{\prime}$ we
find that \[
\alpha_{j^{\prime}}\leq\frac{\Delta-1}{q-\Delta}+\frac{\Delta-1}{(q-\Delta)^{2}}.\]
Then using our assumption that $q\geq2\Delta$ we have \[
\alpha=\max(\alpha_{j},\alpha_{j^{\prime}})\leq\frac{\Delta-1}{q-\Delta}+\frac{\Delta-1}{(q-\Delta)^{2}}\leq\frac{\Delta-1}{\Delta}+\frac{\Delta-1}{\Delta^{2}}=\frac{\Delta^{2}-1}{\Delta^{2}}=1-\frac{1}{\Delta^{2}}<1\]
and we obtain the stated bound on the mixing time by applying Theorem
\ref{thm:main_d}.
\end{proof}

\section{Application: Colouring a Tree\label{sec:tree} }

This section contains the proof of Theorem \ref{thm:tree_block}
which improves the least number of colours required for mixing of
systematic scan on a tree for individual values of $\Delta$. Recall
the definition of the systematic scan $\treeblock$ where the set
of blocks $\Theta$ is defined as follows. Let the block $\Theta_{k}$
contain a site $r$ along with all sites below $r$ in the tree that
are at most $h-1$ edges away from $r$. We call $h$ the \emph{height}
of the blocks and $h$ is defined for each $\Delta$ in Table \ref{tbl:block}.
The set of blocks $\Theta$ covers the sites of the tree and we construct
$\Theta$ such that no block has height less than $h$. $\Pk$ is
the transition matrix for performing a heat-bath move on block $\Theta_{k}$
and hence $\Pk(x)$ is the uniform distribution on the set of configurations
that agree with $x$ off $\Theta_{k}$ and where no edge incident
to a site in $\Theta_{k}$ is monochromatic. The transition matrix
of the Markov chain $\treeblock$ is $\Pi_{k=1}^{m}\Pk$ where $m$
is the number of blocks. 

We will use standard terminology when discussing the structure of
the tree. In particular will say that a site $i$ is a \emph{descendant}
of a site $j$ (or $j$ is a \emph{predecessor} of $i$) if $j$ is
on the simple path from the root of the tree to $i$. We will call
a site $j$ a \emph{child} of a site $i$ (or $i$ is the \emph{parent}
of $j$) if $i$ and $j$ are adjacent and $j$ is a descendant of
$i$. Finally $N_{k}(j)=\{ i\in\partial\Theta_{k}\mid i\textnormal{ is a descendant of }j\}$
is the set of descendants of $j$ on the boundary of $\Theta_{k}$.

Let $(x,y)\in S_{i}$ where $i$ is on the boundary of some block
$\Theta_{k}$. The following lemma will provide upper bounds on the
probability of disagreement at any site in the block. 

\begin{lem}
\label{lemma_block_pr} Let $(x,y)\in S_{i}$ and suppose that $i$
is adjacent to exactly one site in a block $\Theta_{k}$. Then there
exists a coupling $\psi$ of $D_{1}=\Pk(x)$ and $D_{2}=\Pk(y)$ in
which \[
\Prob_{(x^{\prime},y^{\prime})\in\psi}(x_{j}^{\prime}\neq y_{j}^{\prime})\leq\frac{1}{(q-\Delta)^{d(i,j)}}\]
 for all $j\in\Theta_{k}$ where $d(i,j)$ is the edge distance from
$i$ to $j$. 
\end{lem}
\begin{proof}
We construct a coupling $\psi$ of $D_{1}$ and $D_{2}$ based on
the recursive coupling defined in Goldberg et al. \cite{ssm}. The
following definitions are based on Figure \ref{fig:reccouple}. %
\begin{figure}

\caption{The region defined in a boundary pair and the construction of the
subtrees. }

\begin{centering}\label{fig:reccouple}\includegraphics{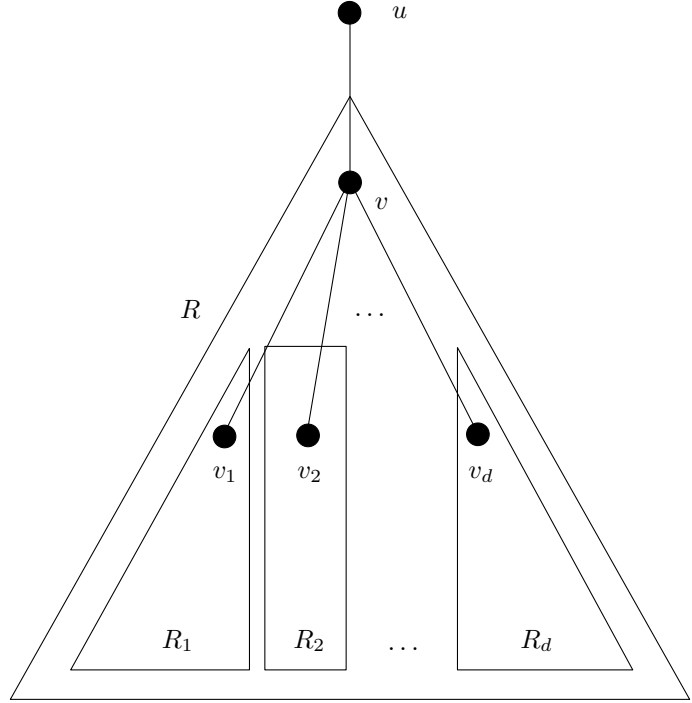}\par\end{centering}
\end{figure}
 Let $R\subseteq V$ be a set of sites. Also let $(X,X^{\prime})$
be a pair of colourings of the sites on the boundary of $R$ (recall
that the boundary of $R$ is the set of sites that are not included
in $R$ but are adjacent to some site in $R$) which use the same
colour for every site, \emph{except} for one site $u$ which is coloured
$l$ in $X$ and $l^{\prime}$ in $X^{\prime}$. We then say that
$A(R,(X,X^{\prime}),u,(l,l^{\prime}))$ is a boundary pair. For a
boundary pair $A(R,(X,X^{\prime}),u,(l,l^{\prime}))$ we let $v\in R$
be the site in $R$ that is adjacent to $u$. We think of $v$ as
the root of $R$ and note that we may need to turn the original tree
``upside down'' in order to achieve this, however the meaning should
be clear. We then label the children (in $R$) of $v$ as $v_{1},\dots,v_{d}$
and let $T=\{ R_{1},\dots,R_{d}\}$ be set the of $d$ subtrees of
$R$ that do not contain site $v$, that is for $R_{k}\in T$ we define
$R_{k}=\{ j\in R\mid j=v_{k}\mbox{ or }j\mbox{ is a descendant of }v_{k}\}$.
Finally let $D$ and $D^{\prime}$ be the uniform distributions on
colourings of $R$ consistent with the boundary colourings $X$ and
$X^{\prime}$ respectively and let $D(v)$ (respectively $D^{\prime}(v)$)
be the uniform distribution on the color at site $v$ induced by $D$
(respectively $D^{\prime}$). Then $\Psi_{R}$ is the recursive coupling
of $D$ and $D^{\prime}$ summarised as follows. 
\begin{enumerate}
\item If $l=l^{\prime}$ then the distributions $D$ and $D^{\prime}$ are
the same and we use the identity coupling, in which the same colouring
is used in both copies. Otherwise we couple $D(v)$ and $D^{\prime}(v)$
greedily to maximise the probability of assigning the same colour
to site $v$ in both distributions. If $R$ consists of a single site
then this completes the coupling.
\item Suppose that the pair of colours $(c,c^{\prime})$ were drawn for
$v$ in the coupling from step 1. For each subtree $R^{\prime}\in\{ R_{1},\dots R_{d}\}$
we have a well defined boundary pair $A(R^{\prime},(X_{R^{\prime}},X_{R^{\prime}}^{\prime}),v,(c,c^{\prime}))$
where $X_{R^{\prime}}$ is the boundary colouring $X$ restricted
to the sites on the boundary of $R^{\prime}$. For each pair of colours
$(c,c^{\prime})$ and $R^{\prime}\in T$ we recursively construct
a coupling $\Psi_{R^{\prime}}(c,c^{\prime})$ of the distributions
induced by the boundary pair $A(R^{\prime},(X_{R^{\prime}},X_{R^{\prime}}^{\prime}),v,(c,c^{\prime}))$. 
\end{enumerate}
Initially we let the boundary pair be $A(R=\Theta_{k},(X=x,Y=y),u=i,(l=x_{i},l^{\prime}=y_{i}))$
and our coupling $\psi$ of $D_{1}$ and $D_{2}$ is thus the recursive
coupling $\Psi_{\Theta_{k}}$ constructed above. 

We prove the statement of the lemma by induction on $d(i,j)$. The
base case is $d(i,j)=1$. Applying Lemma 13 from Goldberg et al. \cite{ssm}
we can upper bound the probability of $x_{j}^{\prime}\neq y_{j}^{\prime}$
where $(x^{\prime},y^{\prime})$ is drawn from $\psi$ by assigning
the worst possible colouring to neighbours of $j$ in $\Theta_{k}$.
Site $j$ has at most $\Delta-1$ neighbours (other than $i$) so
there are at least $q-\Delta$ colours available for $j$ in both
distributions. There is also at most one colour which is valid for
$j$ in $x$ but not in $y$ (and vice versa) so \[
\Prob_{(x^{\prime},y^{\prime})\in\psi}(x_{j}^{\prime}\neq y_{j}^{\prime})\leq\frac{1}{q-\Delta}.\]
Now let $R^{\prime}$ be the subtree of $\Theta_{k}$ containing site
$j$ and let $v$ be the site in $\Theta_{k}$ adjacent to $i$. Assume
that for $d(v,j)=d(i,j)-1$ \[
\Prob_{(x^{\prime},y^{\prime})\in\Psi_{R^{\prime}}(c,c^{\prime})}(x_{j}^{\prime}\neq y_{j}^{\prime})\leq\frac{1}{(q-\Delta)^{d(v,j)}}.\]
Now for $(x,y)\in S_{i}$
\begin{align*}
\Prob_{(x^{\prime},y^{\prime})\in\psi}(x_{j}^{\prime}\neq y_{j}^{\prime})&=\Prob_{(x^{\prime},y^{\prime})\in\Psi_{\Theta_{k}}}(x_{j}^{\prime}\neq y_{j}^{\prime})\\&=\sum_{\begin{subarray}{c} c,c^{\prime} \\
c\neq c^{\prime}\end{subarray}}\Prob_{(x^{\prime},y^{\prime})\in\Psi_{R}}(x_{v}^{\prime}=c,y_{v}^{\prime}=c^{\prime})\Prob_{(x^{\prime},y^{\prime})\in\Psi_{R^{\prime}}(c,c^{\prime})}(x_{j}^{\prime}\neq y_{j}^{\prime}) \\
&\leq\frac{1}{(q-\Delta)^{d(i,j)-1}}\sum_{\begin{subarray}{c} c,c^{\prime} \\
c\neq c^{\prime} \end{subarray}}\Prob_{(x^{\prime},y^{\prime})\in\Psi_{R}}(x_{v}^{\prime}=c,y_{v}^{\prime}=c^{\prime}) \\
&\leq\frac{1}{(q-\Delta)^{d(i,j)}}
\end{align*} where the first inequality is the inductive hypothesis and the last
is a consequence of the base case. 
\end{proof}
We will now use the coupling from Lemma \ref{lemma_block_pr} to define
the coupling $\Psi_{k}(x,y)$ of the distributions $\Pk(x)$ and $\Pk(y)$
for $(x,y)\in S_{i}$. If $i\in\partial\Theta_{k}$ then it is adjacent
to exactly one site in $\Theta_{k}$ and we use the coupling from
Lemma \ref{lemma_block_pr}. If $i\not\in\partial\Theta_{k}$ then
the distributions $\Pk(x)$ and $\Pk(y)$ are the same since we are
using heat-bath updates and so we can use the identity coupling. We
summarise the bounds on $\rho_{i,j}^{k}$ in the following corollary
of Lemma \ref{lemma_block_pr}, which we will use in the proof of
Theorem \ref{thm:tree_block}.

\begin{cor}
\label{cor:tree_block}Let $d(i,j)$ denote the number of edges between
$i$ and $j$. Then for $j\in\Theta_{k}$ \[
\rho_{i,j}^{k}\leq\begin{cases}
\frac{1}{(q-\Delta)^{d(i,j)}} & \textnormal{if }i\in\partial\Theta_{k}\\
0 & \textnormal{otherwise.}\end{cases}\]

\end{cor}
\textbf{\\[1 em]Theorem \ref{thm:tree_block}.} \emph{Let $G$ be
a tree with maximum vertex-degree $\Delta$. If $q\geq f(\Delta)$
where $f(\Delta)$ is specified in Table $\ref{tbl:block}$ for small
$\Delta$ then} \[
\Mix(\treeblock,\epsilon)=O(\log(n\epsilon^{-1})).\]

\begin{proof}
We will use Theorem \ref{thm:main_d} and assign a weight to each
site $i$ such that $w_{i}=\xi^{d_{i}}$ where $d_{i}$ is the edge
distance from $i$ to the root and $\xi$ is defined in Table \ref{tbl:block}
for each $\Delta$. For a block $\Theta_{k}$ and $j\in\Theta_{k}$
we let \[
\alpha_{k,j}=\frac{\sum_{i}w_{i}\rho_{i,j}^{k}}{w_{j}}\]
denote the total weighted influence on site $j$ when updating block
$\Theta_{k}$. For each block $\Theta_{k}$ and each site $j\in\Theta_{k}$
we will upper bound $\alpha_{k,j}$ and hence obtain an upper bound
on $\alpha=\max_{k}\max_{j\in\Theta_{k}}\alpha_{k,j}$. Note from
Corollary \ref{cor:tree_block} that $\rho_{i,j}^{k}=0$ when $i\in\Theta_{k}$
so we only need to bound $\rho_{i,j}^{k}$ for $i\in\partial\Theta_{k}$.

We first consider a block $\Theta_{k}$ that does not contain the
root. The following labels refer to Figure \ref{fig:inf_j_via_v}
in which a solid line is an edge and a dotted line denotes the existence
of a simple path between two sites. Let $p\in\partial\Theta_{k}$
be the predecessor of all sites in $\Theta_{k}$ and $d_{r}-1$ be
the distance from $p$ to the root of the tree i.e., $w_{p}=\xi^{d_{r}-1}$.
The site $r\in\Theta_{k}$ is a child of $p$. Now consider a site
$j\in\Theta_{k}$ which has distance $d$ to $r$, hence $w_{j}=\xi^{d+d_{r}}$
and $d(j,p)=d+1$. From Corollary \ref{cor:tree_block} it then follows
that the weighted influence of $p$ on $j$ when updating $\Theta_{k}$
is at most \[
\rho_{p,j}^{k}\frac{w_{p}}{w_{j}}\leq\frac{1}{(q-\Delta)^{d(j,p)}}\frac{\xi^{d_{r}-1}}{\xi^{d_{r}+d}}=\frac{1}{(q-\Delta)^{d+1}}\frac{1}{\xi^{d+1}}.\]
 Now consider some site $u\in N_{k}(j)$ which is on the boundary
of $\Theta_{k}$. Since $u\in N_{k}(j)$ it has weight $w_{u}=\xi^{d_{r}+h}$
and so $d(j,u)=h-d$. Hence Corollary \ref{cor:tree_block} says that
the weighted influence of $u$ on $j$ is at most \[
\rho_{u,j}^{k}\frac{w_{u}}{w_{j}}\leq\frac{1}{(q-\Delta)^{d(j,u)}}\frac{\xi^{d_{r}+h}}{\xi^{d_{r}+d}}=\frac{1}{(q-\Delta)^{h-d}}\xi^{h-d}.\]
Every site in $\Theta_{k}$ has at most $\Delta-1$ children so the
number of sites in $N_{k}(j)$ is at most $|N_{k}(j)|\leq(\Delta-1)^{h-d}$
and so, summing over all sites $u\in N_{k}(j)$, the total weighted
influence on $j$ from sites in $N_{k}(j)$ when updating $\Theta_{k}$
is at most \[
\sum_{u\in N_{k}(j)}\rho_{u,j}^{k}\frac{w_{u}}{w_{j}}\leq\sum_{u\in N_{k}(j)}\frac{1}{(q-\Delta)^{h-d}}\xi^{h-d}\leq\frac{(\Delta-1)^{h-d}}{(q-\Delta)^{h-d}}\xi^{h-d}.\]

\begin{figure}

\caption{A block in the tree. A solid line indicates an edge and a dotted
line the existence of a path.}

\begin{centering}\label{fig:inf_j_via_v} \includegraphics{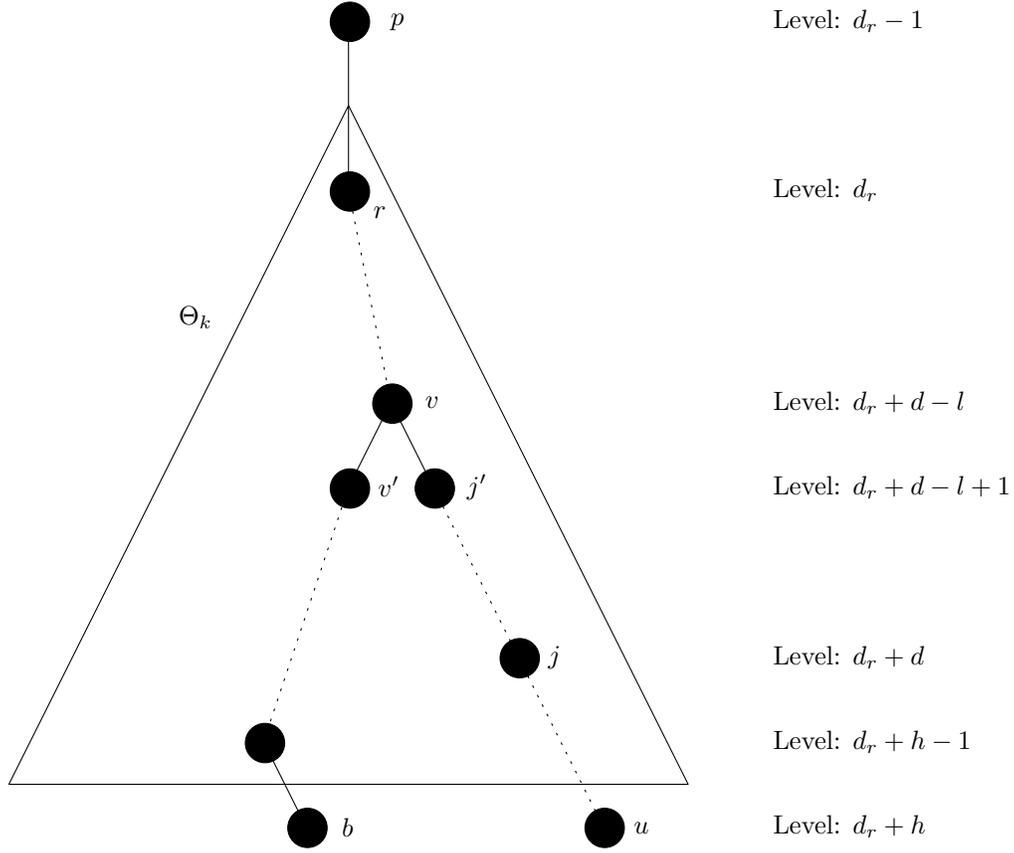} \par\end{centering}
\end{figure}

The influence on $j$ from sites in $\partial\Theta_{k}\setminus\left(N_{k}(j)\cup\{ p\}\right)$
will now be considered. These are the sites on the boundary of $\Theta_{k}$
that are neither descendants or predecessors of $j$. For each site
$v$ between $j$ and $p$, we will bound the influence on site $j$
from sites $b\in N_{k}(v)$ that contain $v$ on the simple path between
$b$ and $j$. We call this the influence on $j$ via $v$. Referring
to Figure \ref{fig:inf_j_via_v} let $v\in\Theta_{k}$ be a predecessor
of $j$ such that $d(j,v)=l$ and observe that $v$ is on level $d_{r}+d-l$
in the tree and also that $1\leq l\leq d$ since $v$ is between $p$
and $j$ in the tree. If $v$ is not the parent of $j$ (that is $l\neq1$)
then let $j^{\prime}$ be the child of $v$ which is also a predecessor
of $j$, that is $j^{\prime}$ is on the simple path from $v$ to
$j$. If $l=1$ we let $j^{\prime}=j$. Also let $v^{\prime}$ be
any child of $v$ other than $j^{\prime}$ and observe that $v^{\prime}$
and $j^{\prime}$ are both on level $d_{r}+d-l+1$. Now let $b\in N_{k}(v^{\prime})$
be a descendant of $v^{\prime}$ and note as before that $w_{b}=\xi^{d_{r}+h}$.
The distance between $b$ and $v^{\prime}$ is \[
d(v^{\prime},b)=d_{r}+h-(d_{r}+d-l+1)=h-d+l-1\]
and so the number of descendants of $v^{\prime}$ is at most $|N_{k}(v^{\prime})|\leq(\Delta-1)^{h-d+l-1}$
since each site has at most $\Delta-1$ children. Site $v$ has at
most $\Delta-2$ children other than $j^{\prime}$ so the number of
sites on the boundary of $\Theta_{k}$ that are descendants of $v$
but not $j^{\prime}$ is at most \[
|N_{k}(v)\setminus N_{k}(j^{\prime})|\leq(\Delta-2)|N_{k}(v^{\prime})|\leq(\Delta-2)(\Delta-1)^{h-d+l-1}.\]
Finally the only simple path from $b$ to $j$ goes via $v$ and the
number of edges on this path is \[
d(j,b)=d(j,v)+d(v,v^{\prime})+d(v^{\prime},b)=l+1+(h-d+l-1)=h-d+2l\]
 so, using Corollary \ref{cor:tree_block}, the weighted influence
of $b$ on site $j$ when updating block $\Theta_{k}$ is at most
\[
\rho_{b,j}^{k}\frac{w_{b}}{w_{j}}\leq\frac{\xi^{d_{r}+h}}{\xi^{d_{r}+d}}\frac{1}{(q-\Delta)^{d(j,b)}}\leq\frac{\xi^{h-d}}{(q-\Delta)^{h-d+2l}}\]
 and summing over all descendants of $v$ (other than descendants
of $j^{\prime}$) on the boundary of $\Theta_{k}$ we find that the
influence on $j$ via site $v$ is at most \begin{equation}
\sum_{b\in N_{k}(v)\setminus N_{k}(j^{\prime})}\rho_{b,j}^{k}\frac{w_{b}}{w_{j}}\leq\sum_{b\in N_{k}(v)\setminus N_{k}(j^{\prime})}\frac{\xi^{h-d}}{(q-\Delta)^{h-d+2l}}\leq\xi^{h-d}\frac{(\Delta-2)(\Delta-1)^{h-d+l-1}}{(q-\Delta)^{h-d+2l}}.\label{eqn:inf_i_via_v}\end{equation}
Summing \eqref{eqn:inf_i_via_v} over $1\leq l\leq d$ gives an upper
bound on the the total weighted influence of sites in $\partial\Theta_{k}\setminus\left(N_{k}(j)\cup\{ p\}\right)$
on site $j$ when updating $\Theta_{k}$ \[
\sum_{b\in\partial\Theta_{k}\setminus\left(N_{k}(j)\cup\{ p\}\right)}\rho_{b,j}^{k}\frac{w_{b}}{w_{j}}\leq\xi^{h-d}\sum_{l=1}^{d}\frac{(\Delta-2)(\Delta-1)^{h-d+l-1}}{(q-\Delta)^{h-d+2l}}\]
 and adding the derived influences we find that the influence on site
$j$ (on level $d_{r}+d$) when updating $\Theta_{k}$ is at most
\begin{eqnarray*}
\alpha_{k,j} & = & \frac{\rho_{p,j}^{k}w_{p}}{w_{j}}+\sum_{u\in N_{k}(j)}\frac{\rho_{u,j}^{k}w_{u}}{w_{j}}+\sum_{b\in\partial\Theta_{k}\setminus\left(N_{k}(j)\cup\{ p\}\right)}\frac{\rho_{b,j}^{k}w_{b}}{w_{j}}\\
 & \leq & \frac{1}{(q-\Delta)^{d+1}}\frac{1}{\xi^{d+1}}+\frac{(\Delta-1)^{h-d}}{(q-\Delta)^{h-d}}\xi^{h-d}+\xi^{h-d}\sum_{l=1}^{d}\frac{(\Delta-2)(\Delta-1)^{h-d+l-1}}{(q-\Delta)^{h-d+2l}}.\end{eqnarray*}

Now consider the block containing the root of the tree, $r$. Let
this be block $\Theta_{0}$ and note that $w_{r}=1$. The only difference
between $\Theta_{0}$ and any other block is that $r$ may have $\Delta$
children. There are at most $\Delta(\Delta-1)^{h-1}$ descendants
of $r$ in $\partial\Theta_{0}$, each of which has weight $\xi^{h}$
so, using Corollary \ref{cor:tree_block}, the weighted influence
on the root is at most \[
\alpha_{0,r}=\sum_{b\in N_{0}(r)}\rho_{b,r}^{0}\frac{w_{b}}{w_{r}}\leq\frac{\Delta(\Delta-1)^{h-1}}{(q-\Delta)^{h}}\xi^{h}.\]

Now consider a site $j$ on level $d\neq0$ in block $\Theta_{0}$.
As in the general case considered above there is an influence of at
most \[
\sum_{b\in N_{0}(j)}\frac{\rho_{b,j}^{0}w_{b}}{w_{j}}\leq\frac{(\Delta-1)^{h-d}}{(q-\Delta)^{h-d}}\xi^{h-d}\]
 on $j$ from the sites in $N_{0}(j)$. Now consider the influence
on site $j$ from $\partial\Theta_{0}\setminus N_{0}(j)$. We first
consider the influence on $j$ via $r$, which is shown in Figure
\ref{fig:infonjroot}.%
\begin{figure}

\caption{The influence on site $j$ via the root. A line denotes an edge and
a dotted line the existence of a simple path.}

\begin{centering}\label{fig:infonjroot}\par\end{centering}

\begin{centering}\includegraphics{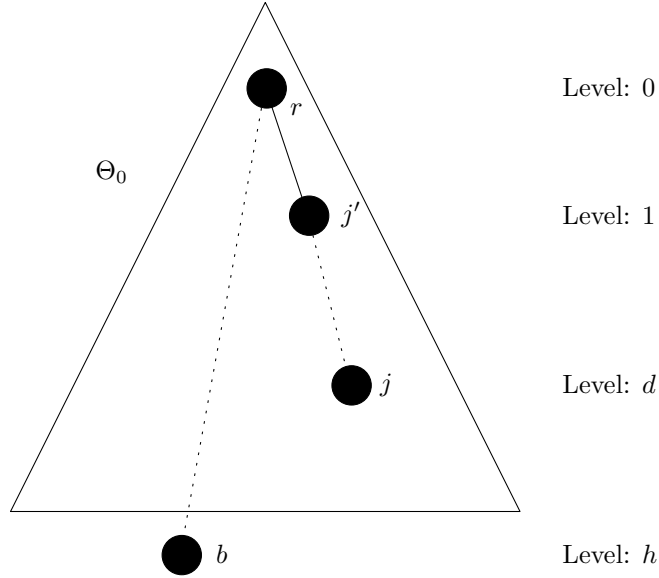}\par\end{centering}
\end{figure}
 Site $r$ has at most $\Delta-1$ children other than the site $j^{\prime}$
which is the child of $r$ that is on the path from $r$ to $j$.
Each child of $r$ has at most $(\Delta-1)^{h-1}$ descendants in
$\partial\Theta_{0}$ and each such descendant has distance $h+d$
to $j$. Hence, from Corollary \ref{cor:tree_block}, the influence
on $j$ via the root is at most \[
\sum_{b\in N_{0}(r)\setminus N_{0}(j^{\prime})}\frac{\rho_{b,j}^{0}w_{b}}{w_{j}}\leq\sum_{b\in N_{0}(r)\setminus N_{0}(j^{\prime})}\frac{\xi^{h}}{\xi^{d}}\frac{1}{(q-\Delta)^{d(b,j)}}\leq\frac{(\Delta-1)^{h}}{(q-\Delta)^{h+d}}\xi^{h-d}.\]

Finally consider then influence on $j$ from the remaining sites,
which are in the set $R=\partial\Theta_{0}\setminus(N_{0}(j)\cup(N_{0}(r)\setminus N_{0}(j^{\prime})))$.
Again consider a site $v\neq r\in\Theta_{0}$ where $v$ is a predecessor
of $j$ and $d(j,v)=l$. In this case we have $1\leq l\leq d-1$ since
$l=d$ is the root which has already been considered. This is the
same situation as arose in the general case considered above (see
Figure \ref{fig:inf_j_via_v}) so \eqref{eqn:inf_i_via_v} is an upper
bound on the influence on $j$ via $v$ and so summing \eqref{eqn:inf_i_via_v}
over $1\leq l\leq d-1$ and adding the other influences on $j$ we
obtain an upper bound on the total weighted influence on site $j$
when updating block $\Theta_{0}$ \begin{eqnarray*}
\alpha_{0,j} & = & \sum_{b\in N_{0}(j)}\frac{\rho_{b,j}^{0}w_{b}}{w_{j}}+\sum_{b\in N_{0}(r)\setminus N_{0}(j^{\prime})}\frac{\rho_{b,j}^{0}w{}_{b}}{w_{j}}+\sum_{b\in R}\frac{\rho_{b,j}^{0}w{}_{b}}{w_{j}}\\
 & \leq & \frac{(\Delta-1)^{h-d}}{(q-\Delta)^{h-d}}\xi^{h-d}+\frac{(\Delta-1)^{h}}{(q-\Delta)^{h+d}}\xi^{h-d}+\xi^{h-d}\sum_{l=1}^{d-1}\frac{(\Delta-2)(\Delta-1)^{h-d+l-1}}{(q-\Delta)^{h-d+2l}}.\end{eqnarray*}

We require $\alpha<1$ which we obtain by satisfying the system of
inequalities given by setting \begin{equation}
\alpha_{k,j}<1\label{eqn:ineqs}\end{equation}
for all blocks $\Theta_{k}$ and sites $j\in\Theta_{k}$. In particular
we need to find an assignment to $\xi$ and $h$ that satisfies \eqref{eqn:ineqs}
given $\Delta$ and $q$. Table~\ref{tbl:block} shows the least
number of colours $f(\Delta)$ required for mixing for small $\Delta$
along with a weight, $\xi$, that satisfies the system of equations
and the required height of the blocks, $h$. These values were verified
by checking the resulting $2h$ inequalities for each $\Delta$ using
Mathematica. The least number of colours required for mixing in the
single site setting is also included in the table for comparison. 
\end{proof}

\section{A Comparison of Influence Parameters\label{sec:comparison-influence}}

We conclude with a discussion of our choice of influence parameter
$\alpha$ denoting the maximum influence on any site in the graph.
As we will be comparing the condition $\alpha<1$ to the corresponding,
but unweighted, conditions in Dyer et al. \cite{dobrushin_scan} and
Weitz \cite{dror_combinatorial} we will let $w_{i}=1$ for each site.
Recall our definitions (letting $w_{i}=1$) of $\rho_{i,j}^{k}$ and
$\alpha$\[
\rho_{i,j}^{k}=\max_{(x,y)\in S_{i}}\{\Prob_{(x^{\prime},y^{\prime})\in\Psi_{k}(x,y)}(x_{j}^{\prime}\neq y_{j}^{\prime})\}\mbox{ and }\alpha=\max_{k}\max_{j\in\Theta_{k}}\sum_{i\in V}\rho_{i,j}^{k}\]
where $\Psi_{k}(x,y)$ is a coupling of the distributions $\Pk(x)$
and $\Pk(y)$. We have previously stated that this is not the standard
way to define the influence of $i$ on $j$ since the coupling is
directly included in the definition of $\rho$. It is worth pointing
out, however, that the corresponding definition in Weitz \cite{dror_combinatorial},
which is also for block dynamics, also makes explicit use of the coupling.
In the single site setting (Dyer et al. \cite{dobrushin_scan}) the
influence of $i$ on $j$, which we will denote $\hat{\rho}_{i,j}$,
is defined by 

\[
\hat{\rho}_{i,j}=\max_{(x,y){}_{i}}\dtv(\mu_{j}(x),\mu_{j}(y))\]
where $\mu_{j}(x)$ is the distribution on spins at site $j$ induced
by $\Pj(x)$. The corresponding condition is $\hat{\alpha}=\max_{j}\sum_{i\in V}\hat{\rho}_{i,j}<1$.
We will show (Lemma \ref{lem:single-rho-same})  that $\hat{\rho}_{i,j}$
is a special case of $\rho_{i,j}^{j}$ when $\Theta_{j}=\{ j\}$ and
$\Psi_{j}(x,y)$ is a coupling minimising the Hamming distance at
site $j$. This will prove our claim that our condition $\alpha<1$
is a generalisation of the single site condition $\hat{\alpha}<1$.
Before establishing this claim we discuss the need to include the
coupling explicitly when working with block dynamics. Consider a pair
of distinct sites $j\in\Theta_{k}$ and $j^{\prime}\in\Theta_{k}$
and a pair of configurations $(x,y)\in S_{i}$. When updating block
$\Theta_{k}$ the dynamics needs to draw a pair of new configurations
$(x^{\prime},y^{\prime})$ from the distributions $\Pk(x)$ and $\Pk(y)$
as previously specified. Hence the interaction between $j$ and $j^{\prime}$
has to be according to these distributions and so it is not possible
to consider the influence of $i$ on $j$ and the influence of $i$
on $j^{\prime}$ separately. In the context of our definition of $\rho$
this means that the influence of $i$ on $j$ and the influence of
$i$ on $j^{\prime}$ have to be defined using the \emph{same coupling}.
This is to say that the coupling $\Psi_{k}(x,y)$ can only depend
on the block $\Theta_{k}$ and the initial pair of configurations
$x$ and $y$, which in turn specify which site is labeled $i$. It
is important to note that the coupling can not depend on $j$, since
otherwise having a small influence on a site would not imply rapid
mixing of systematic scan (or indeed random update). The reason why
we need to make this distinction when working with block dynamics
but not the single site dynamics is that in the single site setting
$\hat{\rho}_{i,j}$ is the influence of site $i$ on $j$ when updating
site $j$ and hence whichever coupling is used must implicitly depend
on $j$. Since the coupling can depend on $j$ in the single site
case it is natural to use the ``optimal'' coupling, which minimises
the probability of having a discrepancy at site $j$. By definition
of total variation distance, the probability of having a discrepancy
at site $j$ under the optimal coupling is $\dtv(\mu_{j}(x),\mu_{j}(y))=\hat{\rho}_{i,j}$
(see e.g. Aldous \cite{aldous_walks}). We will now show that $\hat{\rho}_{i,j}$
is a special case of $\rho_{i,j}^{j}$ in the way described above.

\begin{lem}
\label{lem:single-rho-same}Suppose that for each site $j\in V$ we
have a block $\Theta_{j}=\{ j\}$ and that $\Theta=\{\Theta_{j}\}_{j=1,\dots,n}$.
Also suppose that for each pair $(x,y)\in S_{i}$ of configurations
$\Psi_{j}(x,y)$ is a coupling of $\Pj(x)$ and $\Pj(y)$ in which,
for each $c\in C$, \[
\Prob_{(x^{\prime},y^{\prime})\in\Psi_{j}(x,y)}(x_{j}^{\prime}=y_{j}^{\prime}=c)=\min(\Prob_{\mu_{j}(x)}(c),\Prob_{\mu_{j}(y)}(c))\]
where $\Prob_{\mu_{j}(x)}(c)$ is the probability of drawing colour
$c$ from distribution $\mu_{j}(x)$. Then $\rho_{i,j}^{j}=\hat{\rho}_{i,j}$.
\end{lem}
\begin{proof}
To see that $\Psi_{j}(x,y)$ is a coupling of $\Pj(x)$ and $\Pj(y)$
it is sufficient to observe that $\Prob_{x^{\prime}\in\Pj(x)}(x_{j}^{\prime}=c)=\Prob_{\mu_{j}(x)}(c)$
and similarly $\Prob_{y^{\prime}\in\Pj(y)}(y_{j}^{\prime}=c)=\Prob_{\mu_{j}(y)}(c)$
since $j$ is the only site in $\Theta_{j}$. Thus we have\begin{eqnarray*}
\rho_{i,j}^{j} & = & \max_{(x,y)\in S_{i}}\left\{ \Prob_{(x^{\prime},y^{\prime})\in\Psi_{j}(x,y)}(x_{j}^{\prime}\neq y_{j}^{\prime})\right\} \\
 & = & \max_{(x,y)\in S_{i}}\left\{ 1-\sum_{c\in C}(\Prob_{(x^{\prime},y^{\prime})\in\Psi_{j}(x,y)}(x_{j}^{\prime}=y_{j}^{\prime}=c))\right\} \\
 & = & \max_{(x,y)\in S_{i}}\left\{ 1-\sum_{c\in C}\min(\Prob_{\mu_{j}(x)}(c),\Prob_{\mu_{j}(y)}(c))\right\} \\
 & = & \max_{(x,y)\in S_{i}}\left\{ \sum_{c\in C}\Prob_{\mu_{j}(x)}(c)-\min(\Prob_{\mu_{j}(x)}(c),\Prob_{\mu_{j}(y)}(c))\right\} \\
 & = & \max_{(x,y)\in S_{i}}\left\{ \sum_{c\in C^{+}}\Prob_{\mu_{j}(x)}(c)-\Prob_{\mu_{j}(y)}(c)\right\} \\
 & = & \max_{(x,y)\in S_{i}}\left\{ \frac{1}{2}\sum_{c\in C}|\Prob_{\mu_{j}(x)}(c)-\Prob_{\mu_{j}(y)}(c)|\right\} \\
 & = & \max_{(x,y)\in S_{i}}\dtv(\mu_{j}(x),\mu_{j}(y))\\
 & = & \hat{\rho}_{i,j}\end{eqnarray*}
where $C^{+}=\{ c\mid\Prob_{\mu_{j}(x)}(c)\geq\Prob_{\mu_{j}(y)}(c)\}$.
\end{proof}
Finally we will show that the condition corresponding to $\alpha<1$
in Weitz's paper \cite{dror_combinatorial} does not imply rapid mixing
of systematic scan. Let $B(j)$ be the set of block indices that contain
site $j$ and $b(j)$ the size of this set. Weitz refers to the sum
$\sum_{k\in B(j)}\sum_{i}\rho_{i,j}^{k}$ as the total influence \emph{on}
site $j$ and the parameter representing the maximum influence on
a site, which we denote $\alpha_{W}$ to distinguish it from our own
definition of $\alpha$, is defined as \[
\alpha_{W}=\max_{j}\sum_{k\in B(j)}\sum_{i}\frac{\rho_{i,j}^{k}}{b(j)}.\]
We note that the the single site influence parameter $\hat{\alpha}$
used in Dyer et al. \cite{dobrushin_scan} to prove rapid mixing of
systematic scan is a special case of $\alpha_{W}$ when the coupling
from Lemma \ref{lem:single-rho-same} is used and each site is contained
in exactly one block of size one.

It is proved in Weitz \cite{dror_combinatorial} that the condition
$\alpha_{W}<1$ implies spatial mixing of a random update Markov chain
and hence that the Gibbs measure is unique. We will now show that
the parameters $\alpha$ and $\alpha_{W}$ are different and in particular
that the condition $\alpha_{W}<1$ does not imply rapid mixing of
systematic scan. To show this we exhibit a spin system for which a
systematic scan Markov chain does not mix rapidly but $\alpha_{W}<1$.
It is sufficient to show that a specific systematic scan does not
mix, since Theorem \ref{thm:main_d} states that \emph{any} systematic
scan with a specified set of blocks mixes when $\alpha<1$.

\begin{observation} \label{obs:ring} There exists a spin system
for which $\alpha_{W}<1$ and $\alpha=1$ but systematic scan does
not mix. \end{observation} Consider the following spin system. Let
$G$ be the $n$-vertex cycle and label the sites $0,\dots,n-1$ and
$C$ be the set of $q$ spins. Then $\Theta_{i}$ (which has an associated
transition matrix $\PI$) is the block containing site $i$ and $i+1\mod n$
and it is updated as follows: 

\begin{enumerate}
\item The spin at site $i$ is copied to site $i+1$; 
\item a spin is assigned to site $i$ uniformly at random from the set of
all spins. 
\end{enumerate}
The stationary distribution, $\pi$, of the spin system is the uniform
distribution on all configurations of $G$. Clearly $\PI$ satisfies
property (1) of the update rule, namely that only sites within the
block may change during the update. To see that $\pi$ is invariant
under each $\PI$ observe that site $i+1$ takes the spin of site
$i$ in the original configuration and site $j$ receives a spin drawn
uniformly at random. This ensures that each site has probability $1/q$
of having each spin and that they are independent.

We define the $\rho$ values for this spin system by using the following
coupling. Consider a block $\Theta_{j}$ for update. The spin at site
$j+1$ is deterministic in both copies, and each copy selects the
same colour for site $j$ when drawing uniformly at random from $C$.
First suppose that site $j$ is the discrepancy between two configurations.
Then, since the spin at $j$ is copied to site $j+1$, the spin of
site $j+1$ becomes a disagreement in the coupling and hence $\rho_{j,j+1}^{j}=1$.
The spin at $j$ is drawn uniformly at random from $C$ in both copies
and coupled perfectly so $\rho_{j,j}^{j}=0$. Now suppose that the
two configurations differ at a site $i\neq j$. Then $\rho_{i,j+1}^{j}=0$
since both configurations have the same colour for site $j$, and
$\rho_{i,j}^{j}=0$ since the spins at site $j$ are coupled perfectly.
Using the values of $\rho$ we deduce that \[
\alpha_{W}=\max_{j}\sum_{k\in B(j)}\sum_{i}\frac{\rho_{i,j}^{k}}{b(j)}=\frac{1}{2}\left(\rho_{j-1,j}^{j-1}+\sum_{i\neq j-1}\rho_{i,j}^{j-1}+\sum_{i}\rho_{i,j}^{j}\right)=\frac{1}{2}\]
 and $\alpha=\max_{k}\max_{j\in\Theta_{k}}\sum_{i}\rho_{i,j}^{k}=1$.

Let $\scan$ be the systematic scan Markov chain that updates the
blocks in the order $\Theta_{0},\Theta_{1},\dots,\Theta_{n-1}$. For
each block $\Theta_{i}$ note that if a configuration $y$ is obtained
from updating block $\Theta_{i}$ starting from $x$ then $y_{i+1}=x_{i}$.
Hence when performing the systematic scan, the spin of site $0$ in
the original configuration moves around the ring ending at site $n-1$
before the update of block $\Theta_{n-1}$ moves it on to site 0.
Hence if configuration $x^{\prime}$ is obtained from one complete
scan starting from a configuration $x$ we have $x_{0}^{\prime}=x_{0}$
and the systematic scan Markov chain does not mix since site $0$
will always be assigned the same spin after each complete scan.

\section*{Acknowledgments }

I am grateful to Leslie Goldberg for several useful discussions regarding
technical issues and for providing detailed and helpful comments on
a draft of this article. I would also like to thank Paul Goldberg
for useful comments during the early stages of this work. 

\bibliographystyle{plain}
\bibliography{../references}

\end{document}